\documentclass[a4paper,11pt]{article}
\usepackage[T1]{fontenc}
\usepackage{lmodern}
\usepackage{amsmath,amssymb,amsthm,mathtools}
\usepackage{graphicx,xcolor,array,longtable,booktabs,pdflscape}
\usepackage[shortlabels]{enumitem}
\usepackage[margin=25mm]{geometry}
\usepackage{mathrsfs,microtype}
\usepackage{needspace,eso-pic,afterpage}
\newenvironment{landscapetablepages}{%
  \begin{landscape}%
  \color{black}
  \pagestyle{empty}%
  \AddToShipoutPictureFG{%
    \AtPageLowerLeft{%
      \put(\LenToUnit{\dimexpr\paperwidth-25mm+\footskip\relax},%
           \LenToUnit{.5\paperheight}){%
        \rotatebox{90}{\makebox[0pt][c]{\normalfont\normalsize\color{black}\thepage}}%
      }%
    }%
  }%
}{%
  \clearpage
  \ClearShipoutPictureFG
  \end{landscape}%
}
\usepackage[hidelinks]{hyperref}
\hypersetup{pdfauthor={Himadri Mukherjee, Gunja Sachdeva, Adeep Varghese},pdftitle={Similarity classes of 5x5 matrices over local rings of length two},pdfsubject={Similarity classes in dimension five; reviewed proofs and primary orbit counts}}
\newtheorem{thm}{Theorem}[section]
\newtheorem{prop}[thm]{Proposition}
\newtheorem{lem}[thm]{Lemma}
\newtheorem{cor}[thm]{Corollary}
\theoremstyle{definition}

\theoremstyle{remark}

\DeclareMathOperator{\Stab}{Stab}

\newenvironment{addedblock}{\begingroup}{\endgroup}
\newcommand{\newtext}[1]{{\color{black}#1}}
\newenvironment{blueblock}{\begingroup\color{black}}{\endgroup}
\newcommand{\redtext}[1]{{\color{black}#1}}
\newenvironment{redblock}{\begingroup\color{black}}{\endgroup}

\newcommand{\reviewtext}[1]{{\color{black}#1}}

\newcommand{\edittext}[1]{{\color{black}#1}}
\newenvironment{editblock}{\begingroup\color{black}}{\endgroup}
\newenvironment{reviewblock}{\begingroup\color{black}}{\endgroup}
\newcommand{\proofreadchange}[1]{\textcolor{black}{#1}}
\renewcommand{\arraystretch}{1.2}
\allowdisplaybreaks
\setlist{itemsep=5pt,topsep=4pt}
\numberwithin{equation}{section}
\title{Similarity classes of $5\times5$ matrices\\over local rings of length two}
\author{Himadri Mukherjee \and Gunja Sachdeva \and Adeep Varghese}
\date{}

\newcommand{\C}{\ensuremath{\mathbb{C}}}
\newcommand{\Z}{\ensuremath{\mathbb{Z}}}

\DeclareMathOperator{\Aut}{Aut}

\begin{document}
\begin{addedblock}
\maketitle
\begin{abstract}
\begin{redblock}
\edittext{Let $R$ be a commutative local principal ideal ring of length two whose residue field has arbitrary characteristic. Prasad, Singla and Spallone classified similarity classes in $M_n(R)$ for $n\leq4$ by using module extensions. We extend their classification to $n=5$. We determine explicit extension-orbit representatives and stabilizers for the partitions $(3,2)$, $(2,2,1)$ and $(2,1,1,1)$. We prove that every orbit contains exactly one of the listed representatives. Together with the known primary cases, these results give all similarity classes in $M_5(R)$ and their centralizers. We also determine representatives and stabilizers for the family $(m+1,1,1,1)$ for every $m\geq1$. When the residue field is finite, we obtain the sizes and total number of similarity classes.}
\end{redblock}
\end{abstract}

\medskip
\section{Introduction}\label{sec:intro}
\subsection{The problem and earlier work}
\begin{redblock}
\edittext{For any ring $R$, two matrices $\mathcal A,\mathcal B\in M_n(R)$ are similar if $\mathcal B=g\mathcal A g^{-1}$ for some $g\in\mathrm{GL}_n(R)$. When $R$ is commutative, each matrix $\mathcal A$ defines an $R[t]$-module structure on $R^n$ by $t\cdot v=\mathcal A v$. An $R$-linear isomorphism $v\mapsto gv$ between the modules defined by $\mathcal A$ and $\mathcal B$ is $R[t]$-linear if and only if
\[
 g(\mathcal A v)=\mathcal B(gv)\quad\text{for every }v\in R^n,
\]
or equivalently $g\mathcal A=\mathcal B g$. Thus $\mathcal A$ and $\mathcal B$ are similar if and only if their associated $R[t]$-modules are isomorphic. Over a field, rational canonical form classifies these modules. Over a local ring, reduction modulo the maximal ideal gives only a first step: matrices with the same reduction need not be similar.}

\edittext{We study this problem when $R$ is a commutative local principal ideal ring of length two. Write $(\pi)$ for its maximal ideal and $k=R/(\pi)$, so that $\pi\neq0$ and $\pi^2=0$ in $R$. Examples are $\mathbb Z/\ell^2\mathbb Z$, where $\ell$ is prime, and $k[\varepsilon]/(\varepsilon^2)$. The matrices $0$ and $\pi I_n$ have the same reduction but are not similar. We therefore need to classify the lifts of each similarity class over $k$ and determine their centralizers. If $k$ is finite, the centralizer orders determine the sizes of the similarity classes by the orbit--stabilizer theorem.}

Avni, Onn, Prasad and Vaserstein \cite{AOPV} classify $2\times2$ and $3\times3$ matrices over local principal ideal rings\proofreadchange{. We use their $2\times2$ normal forms below.} Jambor and Plesken \cite{JP} study normal forms over uniserial rings of length two. Prasad, Singla and Spallone \cite{PSS} solve the similarity problem in $M_n(R)$ for $n\leq4$, where $R$ is a commutative local principal ideal ring of length two, by formulating the lifting problem in terms of module extensions. Their correspondence and primary decomposition reduce the problem to explicit group actions indexed by partitions. They also determine the centralizers and list orbit representatives for the families $(m+1,1)$ and $(m+1,1,1)$ \proofreadchange{with $m\geq1$}.

\edittext{Their framework applies to all $n$, but an explicit classification requires solving the orbit problems for all relevant partitions. Even for a fixed $n$, there are several partitions to consider, and finding representatives for each can require a separate, lengthy calculation. As $n$ increases, the number of partitions grows and new orbit problems arise. These calculations are the main difficulty in extending their classification beyond $n=4$.}
\end{redblock}

\subsection{Main results}\label{subsec:main-results}
\begin{reviewblock}
\edittext{In the approach described above, the relevant module extensions are represented by matrices, called extension matrices, on which the corresponding module automorphism group acts. For the partitions $(3,2)$, $(2,2,1)$ and $(2,1,1,1)$, we give orbit representatives and determine their stabilizers. Combined with the known primary cases, these give the following result.}

\begingroup\color{black}
\begin{thm}\label{thm:main}
Let $R$ be a commutative local principal ideal ring of length two.
\begin{enumerate}[label=\textup{\arabic*.},ref=\arabic*,leftmargin=*]
\item\label{main:representatives} The representatives in Theorems~\ref{thm:32}, \ref{thm:221} and \ref{thm:hook}, with $m=1$ in the last theorem, combine with the known primary classifications recalled in Section~\ref{subsec:known-primary} to give a complete set of distinct similarity-class representatives in $M_5(R)$, constructed as in Section~\ref{subsec:completion}. The stabilizers of the extension representatives for these three partitions are given explicitly in Sections~\ref{sec:32}--\ref{sec:hook}.
\item If the residue field $k$ has finite order $Q$, the number of similarity classes in $M_5(R)$ is
\[
 s_5(Q)=Q^{10}+Q^9+3Q^8+4Q^7+8Q^6+8Q^5+10Q^4+6Q^3+6Q^2.
\]
\end{enumerate}
\end{thm}
\endgroup

The stabilizers in part~\ref{main:representatives} determine the centralizers of the corresponding matrices over $R$ through \eqref{eq:centralizer-sequence} and the lifting equation \eqref{eq:lift-centralizer}. \proofreadchange{Section~\ref{sec:enumeration} uses the stabilizer orders to compute the class sizes and the primary orbit counts to compute the number of similarity classes in $M_5(R)$.}
\end{reviewblock}

\begin{redblock}

\edittext{The orbit representatives and stabilizers for $(2,1,1,1)$ are obtained as part of a result for the family $(m+1,1,1,1)$. For any field $K$ and $m\geq1$, Theorem~\ref{thm:hook} gives representatives and stabilizers for the automorphism-group action on the self-extension space of}
\[
 K[\varepsilon]/(\varepsilon^{m+1})\oplus K^3,
\]
\edittext{where $\varepsilon$ acts as zero on $K^3$. If $|K|=q$, the number of orbits is}
\[
 q^m(q^4+2q^3+4q^2+4q+3).
\]

\edittext{The count in Theorem~\ref{thm:main} also has an interpretation over the residue field. By the \proofreadchange{correspondence of Jambor and Plesken \cite{JP}, recalled in \cite[Remark 1.1]{PSS} and \cite{Sharma},} $s_5(Q)$ is the number of simultaneous similarity classes of commuting pairs in $M_5(\mathbb F_Q)$.}
\end{redblock}

\subsection{\edittext{Outline of the paper}}
\begin{editblock}
We use the extension correspondence, primary decomposition and centralizer formulas of Prasad, Singla and Spallone, recalled in Section~\ref{sec:prelim}. Section~\ref{sec:reduction} gives a block reduction used in all three orbit calculations. Sections~\ref{sec:32}--\ref{sec:hook} then give the representatives, compare them, and compute their stabilizers. For $(3,2)$, we first check whether the upper-left entry is invertible and then use the off-diagonal entries to make further reductions. For $(2,2,1)$, we keep track of how changes to the off-diagonal row and column affect the two by two diagonal block. For $(m+1,1,1,1)$, we first put the off-diagonal row and column of the extension matrix into standard forms. We then simplify its three by three diagonal block using transformations that preserve those forms.

Section~\ref{sec:completion} completes the classification in Theorem~\ref{thm:main} by combining the primary cases. It also gives the construction of the matrix representatives over $R$ and describes their centralizers. Section~\ref{sec:enumeration} computes the class sizes and the number of similarity classes in $M_5(R)$.
\end{editblock}

\section{Preliminaries}\label{sec:prelim}
\subsection{Rings and primary modules}
Throughout, $R$ is a commutative local principal ideal ring with maximal ideal $(\pi)$ and
\[
 \reviewtext{\pi\neq0\text{ in }R,\qquad \pi^2=0\text{ in }R},\qquad k=R/(\pi).
\]
Let $\rho:R\longrightarrow k$ be the quotient map, $\rho(r)=r+(\pi)$. We use the same symbol for the map on matrices obtained by reducing each entry modulo $\pi$. A lift of $a\in k$ is an element $\widetilde a\in R$ satisfying $\rho(\widetilde a)=a$. Any two lifts differ by an element of $(\pi)$, so $\pi^2=0$ implies that $\pi\widetilde a$ is independent of the lift. The map $a\mapsto\pi\widetilde a$ identifies the additive group of $k$ with $(\pi)$.

Let $p\in k[t]$ be monic and irreducible, and let $\lambda=(\lambda_1,\ldots,\lambda_r)$ be a partition, with $\lambda_1\geq\cdots\geq\lambda_r>0$. Put
\[
 |\lambda|=\sum_i\lambda_i,\qquad
 S_j=k[t]/(p^j),\qquad K=k[t]/(p),\qquad
 M_\lambda(p)=\bigoplus_i S_{\lambda_i},
\]
and define
\[
 G_\lambda(p)=\Aut_{k[t]}M_\lambda(p),\qquad
 E_\lambda(p)=\operatorname{Ext}^1_{k[t]}(M_\lambda(p),M_\lambda(p)).
\]
\edittext{Replacing $p$ by $cp$, with $c\in k^\times$, gives the same ideals $(p^j)$ and hence the same primary modules. Requiring $p$ to be monic selects one polynomial from these scalar multiples.}

When working with a fixed $p$, we write $M_\lambda$, $G_\lambda$ and $E_\lambda$ in place of $M_\lambda(p)$, $G_\lambda(p)$ and $E_\lambda(p)$, respectively. For a matrix $X$, the notation $p\mid X$ means that every entry is divisible by $p$.

When the residue field $k$ is finite, its order $Q=|k|$ is fixed once $R$ is chosen. The field attached to a primary polynomial $p$ has order
\[
 q=|K|=\lvert k[t]/(p)\rvert=Q^{\deg p}.
\]
For each of the partitions $(3,2)$, $(2,2,1)$ and $(2,1,1,1)$, the primary module $M_\lambda(p)$ has $k$-dimension $5\deg p$. Thus these partitions occur in $M_5(R)$ only when $p$ is linear, in which case $q=Q$.

Each element of $S_j$ has a unique polynomial representative of degree less than $j\deg p$. In Section~\ref{subsec:completion}, we use the coefficients of these polynomials to construct matrices over $R$.

For the calculations in Sections~\ref{sec:32}--\ref{sec:hook}, we use expansions in powers of $p$ with coefficients in $K$. When $p$ is separable, Hensel's lemma gives embeddings $K\hookrightarrow S_j$ such that reduction modulo $p$ sends each embedded element back to itself. We choose these embeddings consistently: the natural map $S_{j+1}\to S_j$ sends each coefficient in the chosen copy of $K$ to the same coefficient in $S_j$. With these identifications, every $u\in S_j$ has a unique expression
\[
 u=u_0+pu_1+\cdots+p^{j-1}u_{j-1},\qquad u_i\in K,
\]
where $u_0$ is the residue of $u$ modulo $p$. We add and multiply the coefficients in $K$ and use $p^j=0$. Equivalently, $S_j\simeq K[\varepsilon]/(\varepsilon^j)$, where $\varepsilon$ is sent to the class of the polynomial $p$ in $S_j$.

The embeddings above preserve multiplication. Simply taking the polynomial representative of degree less than $\deg p$ for each element of $K$ need not do so: multiplication in $K$ is reduced modulo $p$, whereas multiplication in $S_j$ is reduced modulo $p^j$ and may retain nonzero multiples of $p$.

These embeddings are available when $k$ is perfect, in particular when $k$ is finite. For $p=t-a$, the constants in $S_j$ give the required copy of $K=k$ over any residue field. \proofreadchange{In the general hook result, we work over $K[p]$, with $K$ an arbitrary field and $p$ an indeterminate.}

\begin{reviewblock}
We will use cancellation of powers of $p$ only with the corresponding change of modulus. For $1\leq r<j$, multiplication by $p^r$ identifies $S_{j-r}$ with $p^rS_j$ as a $k[t]$-module. Thus
\[
 p^ru=p^rv\text{ in }S_j
 \quad\Longleftrightarrow\quad
 u\equiv v\pmod{p^{j-r}}.
\]
In particular, $pu=pv$ in $S_2$ gives $u_0=v_0$ in $K$, not necessarily $u=v$ in $S_2$. Every subsequent cancellation of $p$ is understood in this sense, entry by entry for matrices.
\end{reviewblock}

\subsection{Lifting and primary decomposition}
For $A\in M_n(k)$, let $M^A$ denote $k^n$ \reviewtext{as a left $k[t]$-module, with $t\cdot v=Av$ for column vectors $v\in k^n$}, and let $G_A=C_{\mathrm{GL}_n(k)}(A)$. The correspondence of \cite[Theorem 2.8]{PSS} gives a bijection
\begin{equation}\label{eq:lifting-bijection}
 \left\{\begin{array}{c}\text{similarity classes in }M_n(R)\\
 \text{whose reductions are similar to }A\end{array}\right\}
 \longleftrightarrow
 G_A\backslash\operatorname{Ext}^1_{k[t]}(M^A,M^A).
\end{equation}
For this bijection, we choose a lift $A_0\in M_n(R)$ of $A$ such that every invertible matrix commuting with $A$ lifts to an invertible matrix commuting with $A_0$. In other words, reduction modulo $\pi$ gives a surjection
\[
 \rho:C_{\mathrm{GL}_n(R)}(A_0)\longrightarrow G_A.
\]
To construct $A_0$, we may take $A$ in primary rational canonical form. For each irreducible polynomial $p$, choose one monic lift $\widetilde p\in R[t]$ and replace every companion block of $p^j$ by the companion matrix of $\widetilde p^{\,j}$, as in \proofreadchange{\cite[Lemma 2.6]{PSS}}. We use this choice of $A_0$ in Section~\ref{subsec:completion}.

The primary decomposition of $M^A$ has the form
\[
 M^A=\bigoplus_p M_{\lambda(p)}(p),\qquad
 \sum_p\deg(p)|\lambda(p)|=n.
\]
\redtext{By \cite[Theorem 4.1]{PSS}, the extension-orbit problem factors over the primary components:}
\begin{equation}\label{eq:primary-orbits}
 G_A\backslash\operatorname{Ext}^1_{k[t]}(M^A,M^A)
 \cong\prod_p G_{\lambda(p)}(p)\backslash E_{\lambda(p)}(p).
\end{equation}
\redtext{\proofreadchange{Only finitely many factors are nontrivial. The stabilizers also factor over these primary components, as in} \cite[Section 4]{PSS}.}

\subsection{Extension matrices and their action}\label{subsec:extension-action}
\redtext{We use the extension-matrix description of \cite[Section 5.1 and Theorem 5.1]{PSS}:}
\[
 \operatorname{Ext}^1_{k[t]}(k[t]/(p^b),k[t]/(p^a))
 \reviewtext{\cong S_a/p^bS_a}\cong k[t]/(p^{\min(a,b)}).
\]
\reviewtext{Here $a,b\geq1$, and these are isomorphisms of $k[t]$-modules.} Thus an element of $E_\lambda$ is represented by an $r\times r$ matrix $E=(E_{ij})$, where
\begin{equation}\label{eq:extension-entries}
 E_{ij}\in S_{\min(\lambda_i,\lambda_j)}.
\end{equation}
An endomorphism $g=(g_{ij})$ of $M_\lambda$ has $g_{ij}$ defined modulo $p^{\lambda_i}$ and divisible by $p^{\max(\lambda_i-\lambda_j,0)}$. For polynomial representatives of its entries, put
\begin{equation}\label{eq:sharp}
 (g^\#)_{ij}\mathrel{\proofreadchange{\equiv}}p^{\lambda_j-\lambda_i}g_{ij}
 \proofreadchange{\pmod{p^{\lambda_j}}}.
\end{equation}
The divisibility condition makes this expression well-defined. By \cite[Theorem 5.2]{PSS}, the action of $G_\lambda$ on $E_\lambda$ is given by
\begin{equation}\label{eq:action}
 g\cdot E=gE(g^\#)^{-1},
\end{equation}
\proofreadchange{where the $(i,j)$ entry belongs to $S_{\min(\lambda_i,\lambda_j)}$ as in \eqref{eq:extension-entries}.} Equivalently,
\begin{equation}\label{eq:intertwine}
 E'=g\cdot E\quad\Longleftrightarrow\quad gE=E'g^\#.
\end{equation}
\redtext{When all parts of $\lambda$ are equal, the action is ordinary similarity, as in \cite[Theorem 8.3]{PSS}.}

We write
\[
 \Stab_{G_\lambda}(E)=\{g\in G_\lambda:gE=Eg^\#\},\qquad
 \mathcal O_E=G_\lambda\cdot E.
\]
\edittext{When $|K|=q$, write $c_\lambda(q)$ for the number of $G_\lambda$-orbits on $E_\lambda$, and put $c_\varnothing(q)=1$. In the counting tables, $N$ is the number of representatives in a family and $H$ is their common stabilizer order. The columns $|G_\lambda|/H$ and $N|G_\lambda|/H$ give the orbit size and the number of elements in that family of orbits. Summing the last column gives $|E_\lambda|$.}

\subsection{Stabilizers and centralizers}
\redtext{For the primary decomposition $\mathcal A=\bigoplus_p\mathcal A_p$ of \cite[Section 4 and Remark 4.2]{PSS}, with distinct indexing polynomials and $\mathcal A_p\in M_{n_p}(R)$, the centralizer decomposes as}
\begin{equation}\label{eq:centralizer-primary-product}
 C_{\mathrm{GL}_n(R)}(\mathcal A)
 \cong\prod_p C_{\mathrm{GL}_{n_p}(R)}(\mathcal A_p),
 \qquad n=\sum_p n_p.
\end{equation}

\redtext{Let $\mathcal A\in M_n(R)$ correspond to an extension class $E$ above $A=\rho(\mathcal A)$. Equation (3.1) and Theorem 3.1 of \cite{PSS} give the exact sequence}
\begin{equation}\label{eq:centralizer-sequence}
 1\longrightarrow(C_{M_n(k)}(A),+)
 \xrightarrow{\,X\mapsto I+\pi\widetilde X\,}
 C_{\mathrm{GL}_n(R)}(\mathcal A)
 \xrightarrow{\rho}\Stab_{G_A}(E)\longrightarrow1.
\end{equation}
If $|k|=Q$, put
\begin{equation}\label{eq:delta}
 \delta(A)=\sum_p\deg(p)\sum_{i,j}\min(\lambda_i(p),\lambda_j(p)).
\end{equation}
\redtext{The order formula of \cite[Corollary 3.2 and Theorem 9.5]{PSS}, combined with primary decomposition, is}
\begin{equation}\label{eq:global-centralizer}
 |C_{\mathrm{GL}_n(R)}(\mathcal A)|
 =Q^{\delta(A)}\prod_p|\Stab_{G_{\lambda(p)}(p)}(E_p)|.
\end{equation}
\proofreadchange{Here $E_p$ is the $p$-primary component of $E$.}

\edittext{Section~\ref{subsec:completion} gives linear equations for the lifts of each stabilizer element to the centralizer.}

\section{A block reduction}\label{sec:reduction}
\end{addedblock}

\begin{addedblock}
Let $\lambda=(m^r,n^s)$, with $m>n\geq1$ and $\ell=m-n$. Then
\[
 G_\lambda=\left\{\begin{pmatrix}A&p^\ell B\\C&D\end{pmatrix}:\begin{array}{l}
 A\in\mathrm{GL}_r(S_m),\quad D\in\mathrm{GL}_s(S_n),\\
 \newtext{B\in M_{r\times s}(S_n),\quad C\in M_{s\times r}(S_n)}
 \end{array}\right\},
\]
and
\[
 E_\lambda=\left\{\begin{pmatrix}X&Y\\Z&W\end{pmatrix}:\begin{array}{l}
 X\in M_r(S_m),\quad W\in M_s(S_n),\\
 \newtext{Y\in M_{r\times s}(S_n),\quad Z\in M_{s\times r}(S_n)}
 \end{array}\right\}.
\]
\end{addedblock}
\begin{addedblock}
\Needspace{7\baselineskip}
\begin{blueblock}
\edittext{The next proposition is used in the invertible-block case of each classification.}
\end{blueblock}

\begin{prop}\label{1stt}
Let $\lambda=(m^r,n^s)$, where $m>n\geq1$ and $\ell=m-n$. If
$E=\left(\begin{smallmatrix}X&Y\\Z&W\end{smallmatrix}\right)\in E_\lambda$
and $X\in\mathrm{GL}_r(S_m)$, then $E$ lies in the same $G_\lambda$-orbit as a block-diagonal matrix.
\end{prop}
\begin{proof}
We look for $g=\left(\begin{smallmatrix}A&p^\ell B\\C&D\end{smallmatrix}\right)$ and $E'=\operatorname{diag}(X',W')$ such that $gE=E'g^\#$. The four block equations are
\begin{align*}
 AX+p^\ell BZ&=X'A,& CY+DW&=W'D,\\
 AY+p^\ell BW&=X'B,& CX+DZ&=p^\ell W'C.
\end{align*}
Take $A=I_r$ and $D=I_s$. Let $X_n$ be the reduction of $X$ to $M_r(S_n)$. The diagonal equations give $X'=X+p^\ell BZ$ and $W'=W+CY$. Reducing the off-diagonal equations modulo $p^\ell$ gives
\[
 B\equiv X_n^{-1}Y,\qquad C\equiv\newtext{-ZX_n^{-1}}\pmod{p^\ell}.
\]
\begin{blueblock}
These are the first approximations. To obtain the required equalities in $S_n$, use the full equations
\begin{align*}
 B&=X_n^{-1}(Y+p^\ell BW-p^\ell BZB),\\
 C&=(-Z+p^\ell WC+p^\ell CYC)X_n^{-1}.
\end{align*}
Starting from $B=C=0$, repeatedly substitute into the right-hand sides. If two inputs agree modulo $p^j$, their next values agree modulo $p^{j+\ell}$. Since $p^n=0$ in $S_n$, successive iterates eventually agree in $S_n$ and hence solve both equations exactly. The resulting matrix $g$ is invertible because its Schur complement $I_s-p^\ell CB$ is congruent to $I_s$ modulo $p$. Thus it gives the required block-diagonal representative.
\end{blueblock}
\end{proof}
\end{addedblock}

\section{The partition \texorpdfstring{$(3,2)$}{(3,2)}}\label{sec:32}
In this section, $M_\lambda=S_3\oplus S_2$. \newtext{A subscript $0$ denotes reduction modulo $p$.} Thus
\begin{align*}
 E_\lambda&=\left\{\begin{pmatrix}x&y\\z&w\end{pmatrix}:x\in S_3,\ y,z,w\in S_2\right\},\\
 G_\lambda&=\left\{\begin{pmatrix}a&pb\\c&d\end{pmatrix}:a\in S_3^\times,\ d\in S_2^\times,\ b,c\in S_2\right\}.
\end{align*}
\Needspace{8\baselineskip}
\edittext{Write $X=\left(\begin{smallmatrix}x&y\\z&w\end{smallmatrix}\right)$ and $X'=\left(\begin{smallmatrix}x'&y'\\z'&w'\end{smallmatrix}\right)$. The matching equation $gX=X'g^\#$ is equivalent to}
\begin{align}
 ax+pbz&\mathrel{\proofreadchange{\equiv}}ax'+pcy'\pmod{p^3},\label{eq:32a}\\
 ay+pbw&\mathrel{\proofreadchange{\equiv}}bx'+dy'\pmod{p^2},\label{eq:32b}\\
 cx+dz&\mathrel{\proofreadchange{\equiv}}az'+pcw'\pmod{p^2},\label{eq:32c}\\
 cy+dw&\mathrel{\proofreadchange{\equiv}}bz'+dw'\pmod{p^2}.\label{eq:32d}
\end{align}

\begin{thm}\label{thm:32}
The following families form a complete set of distinct orbit representatives:
\begin{enumerate}[label=(\roman*)]
\item $\operatorname{diag}(x,w)$, where $x\in S_3^\times$ and $w\in S_2$.\label{o321}
\item $\operatorname{diag}(p\newtext{\xi},w)$, where $\newtext{\xi},w\in S_2$.\label{o321sing}
\item $\begin{pmatrix}0&1\\z&w\end{pmatrix}$, where $z,w\in S_2$.\label{o322}
\item $\begin{pmatrix}0&py\\1&w\end{pmatrix}$, where $y\in K$ and $w\in S_2$.\label{o323}
\item $\begin{pmatrix}px&p\\pz&x+pw\end{pmatrix}$, where $x,z,w\in K$.\label{o324}
\item $\begin{pmatrix}px&0\\p&x+pw\end{pmatrix}$, where $x,w\in K$.\label{o325}
\end{enumerate}
\end{thm}

\begin{proof}
\newtext{\textbf{Exhaustiveness of the list.}}
Let $\newtext{X}=\left(\begin{smallmatrix}x&y\\z&w\end{smallmatrix}\right)\in E_\lambda$.

\textbf{Case 1: $p\nmid x$.}
By Proposition~\ref{1stt}, $\newtext{X}$ can be reduced to a diagonal matrix. Its upper-left entry remains a unit, so we obtain family~\ref{o321}.

\textbf{Case 2: $p\mid x$.}
Write $x=p\newtext{\xi}$, where $\newtext{\xi}\in S_2$. \newtext{We write $\xi=\xi_0+p\xi_1$, with $\xi_0,\xi_1\in K$ using the fixed coefficient representatives; in particular, \proofreadchange{$\xi_0\equiv\xi\pmod{p}$}.} We use the following subcases.

\textbf{Case 2.1: $p\nmid y$.}
The matrix
\[
 g=\begin{pmatrix}1&0\\\newtext{\xi}&y\end{pmatrix}
\]
takes $\newtext{X}$ to a representative of family~\ref{o322}.
\newtext{More explicitly, substituting in \eqref{eq:32a}--\eqref{eq:32d} gives}
\[
 \newtext{g\cdot \newtext{X}=\begin{pmatrix}0&1\\yz-p\newtext{\xi} w&w+\newtext{\xi}\end{pmatrix}.}
\]
\proofreadchange{Here and below, the upper-left entry belongs to $S_3$ and the other entries belong to $S_2$.}

\textbf{Case 2.2: $p\mid y$ and $p\nmid z$.}
Write $y=py_1$, with $y_1\in K$. The matrix
\[
 g=\begin{pmatrix}1&-p\newtext{\xi}/z\\0&1/z\end{pmatrix}
\]
gives family~\ref{o323}. \newtext{Indeed,}
\[
 \newtext{g\cdot \newtext{X}=\begin{pmatrix}0&p(zy_1-\newtext{\xi} w)\\1&w+\newtext{\xi}\end{pmatrix}.}
\]
\newtext{Only the residue of $zy_1-\newtext{\xi} w$ occurs in the upper-right entry. All divisions in this step are by the unit $z$.}

\textbf{Case 2.3: $p\mid y,z$, $p\nmid w$, and $p\nmid\proofreadchange{(\xi-w)}$.}
Write $y=py_1$ and $z=pz_1$, with $y_1,z_1\in K$, and put $\Delta=\newtext{\xi}-w$. The matrix
\[
 g=\begin{pmatrix}
 1&py_1/\Delta\\
 -z_1/\Delta&1-py_1z_1/\Delta^2
 \end{pmatrix}
\]
reduces $\newtext{X}$ to a diagonal matrix, giving family~\ref{o321sing}.
\edittext{The resulting diagonal entries are}
\[
 \newtext{x'=p\newtext{\xi}+\frac{p^2y_1z_1}{\Delta}\quad\text{in }S_3,
 \qquad w'=w-\frac{py_1z_1}{\Delta}\quad\text{in }S_2.}
\]
\newtext{These entries and $y'=z'=0$ satisfy all four matching equations. The denominators are units, and the diagonal entries of $g$ have nonzero residues, so $g\in G_\lambda$.}

\textbf{Case 2.4: $p\mid y,z$, $p\nmid w$, and $p\mid\proofreadchange{(\xi-w)}$.}
\newtext{Here $\xi_0=w_0\neq0$. Write}
\[
 \newtext{x=p\xi_0+p^2\xi_1},\qquad w=\newtext{\xi_0}+pw_1,\qquad
 y=py_1,\quad z=pz_1,
\]
where all five coefficients are in $K$ and $\newtext{\xi_0}\neq0$.
If $y_1\neq0$, use
\[
 g=\begin{pmatrix}y_1^{-1}&0\\\newtext{\xi_1}y_1^{-1}&1\end{pmatrix}.
\]
\newtext{Its image is}
\[
 \newtext{\begin{pmatrix}p\newtext{\xi_0}&p\\py_1z_1&\newtext{\xi_0}+p(w_1+\newtext{\xi_1})\end{pmatrix},}
\]
which belongs to family~\ref{o324}. If $y_1=0$ and $z_1\neq0$, use
\[
 g=\begin{pmatrix}z_1&\newtext{-p\newtext{\xi_1}}\\0&1\end{pmatrix}.
\]
\newtext{Substitution gives the image}
\[
 \newtext{\begin{pmatrix}p\newtext{\xi_0}&0\\p&\newtext{\xi_0}+p(w_1+\newtext{\xi_1})\end{pmatrix},}
\]
which belongs to family~\ref{o325}. If $y_1=z_1=0$, $\newtext{X}$ is already diagonal.

\textbf{Case 2.5: $p\mid \newtext{X}$.}
Write
\[
 \newtext{X}=p\newtext{X_1},\qquad
 \newtext{X_1}=\begin{pmatrix}\newtext{\xi}&y_1\\z_1&w_1\end{pmatrix},
 \quad \newtext{\xi}\in S_2,\quad y_1,z_1,w_1\in K.
\]
\newtext{\reviewtext{Cancelling $p$ in the matching equations changes the top-left modulus from $p^3$ to $p^2$ and the other moduli from $p^2$ to $p$. The resulting equations are precisely those for the partition $(2,1)$.} Every automorphism for this shorter partition lifts to one in $G_{(3,2)}$. If $\newtext{\xi}$ is a unit, Proposition~\ref{1stt} makes $\newtext{X_1}$ diagonal. Otherwise write $\newtext{\xi}=p\newtext{\xi_1}$. If $y_1\neq0$, the transformation of Case~2.1 for $(2,1)$ gives a matrix $\left(\begin{smallmatrix}0&1\\z'&w'\end{smallmatrix}\right)$. If $y_1=0$ and $z_1\neq0$, Case~2.2 gives $\left(\begin{smallmatrix}0&0\\1&w'\end{smallmatrix}\right)$. If both vanish, the matrix is diagonal. Multiplication by $p$ gives family~\ref{o321sing}, or families~\ref{o324} and \ref{o325} with $x=0$.}

\newtext{\textbf{Comparison of representatives.}}
Suppose $gX=X'g^\#$ with $g=\left(\begin{smallmatrix}a&pb\\c&d\end{smallmatrix}\right)\in G_\lambda$. In particular, $a_0,d_0\neq0$. Equation~\eqref{eq:32a} modulo $p$ shows that the residue of the upper-left entry is unchanged. Thus family~\ref{o321} cannot meet any other family.

If a matrix of family~\ref{o321sing} is in the same orbit as one of family~\ref{o322}, equation~\eqref{eq:32b} gives $d_0=0$, a contradiction. Comparing families~\ref{o321sing} and \ref{o323}, equation~\eqref{eq:32c} instead gives $a_0=0$. For families~\ref{o321sing} and \ref{o324}, write $X=\operatorname{diag}(p\newtext{\xi},w)$ and let $x'$ be the upper-left parameter in the second representative. Equation~\eqref{eq:32d} gives $w_0=x'$, whereas \reviewtext{\eqref{eq:32b} gives $p\bigl(b(w-x')-d\bigr)=0$ in $S_2$. By the cancellation rule in Section~\ref{sec:prelim}, the expression in parentheses has zero residue, so}
\[
 b_0(w_0-x')=d_0,
\]
again a contradiction. For families~\ref{o321sing} and \ref{o325}, equation~\eqref{eq:32a} gives $\newtext{\xi}=x'$ in $S_2$. Substitution in \eqref{eq:32c} then gives $pa=0$ in $S_2$, which contradicts $a_0\neq0$.

Comparing family~\ref{o322} with any of \ref{o323}, \ref{o324} or \ref{o325}, equation~\eqref{eq:32b} modulo $p$ gives $a_0=0$. Comparing family~\ref{o323} with \ref{o324} or \ref{o325}, equation~\eqref{eq:32c} modulo $p$ gives $d_0=0$. Finally, for families~\ref{o324} and \ref{o325}, equation~\eqref{eq:32d} modulo $p$ gives $x=x'$, since these parameters lie in $K$. \reviewtext{Cancelling $p$ in \eqref{eq:32b} and passing from $S_2$ to $K$ now gives $a_0=0$.} Hence no two different families meet the same orbit.

It remains to compare two representatives within the same family. For the diagonal families~\ref{o321} and \ref{o321sing}, equations~\eqref{eq:32a} and \eqref{eq:32d}, together with the invertibility of $a$ and $d$, give equality of the diagonal entries.

For family~\ref{o322}, equations~\eqref{eq:32a} and \eqref{eq:32b} give
\[
 c=bz,\qquad d=a+pbw\qquad\text{in }S_2.
\]
\begin{blueblock}
Substituting these identities in \eqref{eq:32c} and \eqref{eq:32d} gives the two equations
\[
 \begin{pmatrix}a&pbz\\b&d\end{pmatrix}
 \binom{z-z'}{w-w'}=\binom00\qquad\text{in }S_2.
\]
Here $a$ denotes its reduction to $S_2$. The determinant $ad-pb^2z$ has nonzero residue $a_0d_0$, so this matrix is invertible. Therefore $z=z'$ and $w=w'$.
\end{blueblock}

For family~\ref{o323}, equation~\eqref{eq:32a} gives $b=py'c$ in $S_2$, and hence $b_0=\newtext{0}$. Equation~\eqref{eq:32c} gives $d_0=a_0$. \reviewtext{Cancelling $p$ in \eqref{eq:32b} and passing from $S_2$ to $K$ therefore gives} $a_0y=d_0y'$, so $y=y'$. Now \eqref{eq:32d} becomes
\[
 pcy+dw=pyc+dw',
\]
and the invertibility of $d$ gives $w=w'$.

For family~\ref{o324}, equation~\eqref{eq:32d} modulo $p$ gives $x=x'$. Substitution in \eqref{eq:32b}, \reviewtext{followed by cancellation of $p$ and reduction to $K$, gives $a_0=d_0$}. Equation~\eqref{eq:32c} similarly gives $d_0z=a_0z'$, so $z=z'$. Equation~\eqref{eq:32a} now gives \newtext{$c_0=b_0z$ in $K$}. Using this in \eqref{eq:32d}, \reviewtext{after subtracting the common term $dx$ and cancelling $p$ to pass from $S_2$ to $K$}, gives
\[
 c_0+d_0w=b_0z+d_0w',
\]
so $w=w'$.

\begin{blueblock}
For family~\ref{o325}, equation~\eqref{eq:32d} modulo $p$ first gives $x=x'$. Equation~\eqref{eq:32a} then gives $p^2b=0$ in $S_3$, hence $b_0=0$. \reviewtext{Cancelling $p$ in \eqref{eq:32c} gives $d_0=a_0$ in $K$.} Finally, \eqref{eq:32d} gives $d_0(w-w')=b_0=0$, so $w=w'$. Together with exhaustiveness, this proves that each orbit contains exactly one listed matrix.
\end{blueblock}
\end{proof}

\subsection{Stabilizers and counts}

\begin{blueblock}
\edittext{For each representative $X$ in Theorem~\ref{thm:32}, the table below gives $\Stab_{G_\lambda}(X)$. In every row, $a\in S_3^\times$. Write}
\[
 a=a_0+pa_1+p^2a_2,\qquad
 \proofreadchange{a^{(2)}=a_0+pa_1\in S_2,\qquad a^{(2)}\equiv a\pmod{p^2}},
\]
\edittext{and \proofreadchange{let $b_0\in K$ satisfy $b_0\equiv b\pmod{p}$.} The free parameters in each row vary independently. The parameters $\xi,y,z,w$ belong to the fixed representative.}
\end{blueblock}
\begin{center}\small
\renewcommand{\arraystretch}{1.65}
\begin{tabular}{c c l}
\toprule
Type & Matrices in $\Stab_{G_\lambda}(X)$ & Free parameters besides $a\in S_3^\times$\\\midrule
$\text{(i)}$ & $\begin{pmatrix}a&0\\0&d\end{pmatrix}$ & $d\in S_2^\times$\\[8pt]
$\text{(ii)},\ \xi_0=w_0$ & $\begin{pmatrix}a&pb\\c&d\end{pmatrix}$ & $b,c\in S_2,\ d\in S_2^\times$\\[8pt]
$\text{(ii)},\ \xi_0\neq w_0$ & $\begin{pmatrix}a&p^2b_1\\pc_1&d\end{pmatrix}$ & $b_1,c_1\in K,\ d\in S_2^\times$\\[8pt]
$\text{(iii)}$ & $\begin{pmatrix}a&pb\\bz&\newtext{a^{(2)}}+pbw\end{pmatrix}$ & $b\in S_2$\\[8pt]
$\text{(iv)}$ & $\begin{pmatrix}a&p^2cy\\c&\newtext{a^{(2)}}+pcw\end{pmatrix}$ & $c\in S_2$\\[8pt]
$\text{(v)}$ & $\begin{pmatrix}a&pb\\\newtext{b_0z+pc_1}&\newtext{a_0+pd_1}\end{pmatrix}$ & $b\in S_2,\ c_1,d_1\in K$\\[8pt]
$\text{(vi)}$ & $\begin{pmatrix}a&p^2b_1\\c&\newtext{a_0+pd_1}\end{pmatrix}$ & $c\in S_2,\ b_1,d_1\in K$\\[8pt]\bottomrule
\end{tabular}
\end{center}
\begin{blueblock}
\proofreadchange{The upper-right entries belong to $pS_3$. The lower-left and lower-right entries belong to $S_2$.} In types~(iii)--(vi), the lower-right entry has residue $a_0\neq0$, so no further invertibility condition is needed. In particular, $a_2$ is free in every row. In types~(v) and (vi), $d_1$ is also independent of $a_1$: only the residues of the two diagonal entries must coincide.

\edittext{Set $X'=X$ in \eqref{eq:32a}--\eqref{eq:32d}. In type (i), the unit $x$ forces $b=c=0$. In type (ii), the conditions are $pb(\xi-w)=pc(\xi-w)=0$ in $S_2$. They impose no restriction if $\xi_0=w_0$, and otherwise give $b,c\in pS_2$. For types (iii)--(vi), the coefficient relations obtained in the comparison proof give the displayed matrices. Substitution in the four matching equations shows that these conditions are also sufficient.}

\end{blueblock}
\edittext{For $|K|=q$, we have $|G_\lambda|=q^7(q-1)^2$. Using the table conventions of Section~\ref{subsec:extension-action}, we obtain}
\begin{center}\small
\begin{tabular}{c c c c c}
\toprule
Type & $N$ & $H$ & $|\mathcal O|=|G_\lambda|/H$ & $N|\mathcal O|$\\\midrule
$\text{(i)}$ & $q^4(q-1)$ & $q^3(q-1)^2$ & $q^4$ & $q^8(q-1)$\\
$\text{(ii)},\ \xi_0=w_0$ & $q^3$ & $q^7(q-1)^2$ & $1$ & $q^3$\\
$\text{(ii)},\ \xi_0\neq w_0$ & $q^3(q-1)$ & $q^5(q-1)^2$ & $q^2$ & $q^5(q-1)$\\
$\text{(iii)}$ & $q^4$ & $q^4(q-1)$ & $q^3(q-1)$ & $q^7(q-1)$\\
$\text{(iv)}$ & $q^3$ & $q^4(q-1)$ & $q^3(q-1)$ & $q^6(q-1)$\\
$\text{(v)}$ & $q^3$ & $q^6(q-1)$ & $q(q-1)$ & $q^4(q-1)$\\
$\text{(vi)}$ & $q^2$ & $q^6(q-1)$ & $q(q-1)$ & $q^3(q-1)$\\\bottomrule
\end{tabular}
\end{center}

\begin{blueblock}
The orbit sizes sum, with these multiplicities, to
\[
 (q^9-q^8)+q^3+(q^6-q^5)+(q^8-q^7)+(q^7-q^6)
 +(q^5-q^4)+(q^4-q^3)=q^9.
\]
\edittext{Summing the column $N$ gives}
\begin{equation}\label{eq:c32}
 c_{(3,2)}(q)=q^5+q^4+2q^3+q^2.
\end{equation}

\end{blueblock}
\section{The partition \texorpdfstring{$(2,2,1)$}{(2,2,1)}}\label{sec:221}
Here $M_\lambda=S_2^2\oplus K$. We write
\[
 E=\begin{pmatrix}X&\bar{v}\\\bar{u}^{\,t}&w\end{pmatrix},\qquad
 g=\begin{pmatrix}A&p\bar{b}\\\bar{c}^{\,t}&e\end{pmatrix},
\]
where $X\in M_2(S_2)$, $A\in\mathrm{GL}_2(S_2)$, $\bar u,\bar v,\bar b,\bar c\in K^2$, $w\in K$ and $e\in K^\times$. \newtext{The bars denote vectors. We write $\bar v=(v_1,v_2)^t$ and $\bar u=(u_1,u_2)^t$, so the lower-left block is the row $\bar u^t$.} Put
\[
 C(a,b)=\begin{pmatrix}0&1\\a&b\end{pmatrix},\qquad
 e_1=(1,0)^t,\quad e_2=(0,1)^t.
\]
\newtext{Block-diagonal changes of basis act on $X$ by ordinary similarity, so we first recall its normal forms. The off-diagonal blocks can change $X$ further, as the matching equations below show.}

\Needspace{12\baselineskip}
\begin{thm}[Avni--Onn--Prasad--Vaserstein {\cite[Theorem~2.2]{AOPV}}]\label{nf}
Every matrix in $M_2(S_2)$ is similar to exactly one matrix in the following list:
\begin{enumerate}[label=(\roman*)]
\item $dI_2$, with $d\in S_2$;
\item $C(a,b)$, with $a,b\in S_2$;
\item $dI_2+pC(a,b)$, with $d,a,b\in K$.
\end{enumerate}
\end{thm}
\begin{blueblock}
\edittext{These are the cases $j=2,0,1$, respectively, of the cited theorem, using the coefficient representatives fixed in Section~\ref{sec:prelim}.}
\end{blueblock}

The matching equation $gE=E'g^\#$ gives
\begin{align}
 AX+p\bar{b} \bar{u}^{\,t}&=X'A+p\newtext{\bar{v}'}\bar{c}^{\,t},\label{me1}\\
 A\bar{v}&=X'\bar{b}+e\newtext{\bar{v}'},\label{me2}\\
 \bar{c}^{\,t} X+e\bar{u}^{\,t}&=\newtext{\bar{u}'^{\,t}}A,\label{me3}\\
 \bar{c}^{\,t} \bar{v}+ew&=\newtext{\bar{u}'^{\,t}}\bar{b}+ew'.\label{me4}
\end{align}
The first equation is in $M_2(S_2)$ and the remaining equations are over $K$. \newtext{Reduction of \eqref{me1} modulo $p$ shows that the reductions of $X$ and $X'$ are similar. We therefore consider the following three cases separately.}

\subsection{\texorpdfstring{Case 1: $X$ is invertible}{Case 1: X is invertible}}
\begin{prop}
If $X$ is invertible, then $E$ can be reduced to $\operatorname{diag}(X',w')$, where $X'$ is one of the invertible normal forms in \newtext{Theorem~\ref{nf}}.
\end{prop}
\begin{proof}
Proposition~\ref{1stt} removes the off-diagonal blocks. A block-diagonal change of basis then puts $X'$ into one of the forms in \newtext{Theorem~\ref{nf}}. \newtext{Its reduction is still invertible. These representatives will be types 12--14 in Theorem~\ref{thm:221}.}
\end{proof}

\subsection{\texorpdfstring{Case 2: $X$ is noninvertible and $p\nmid X$}{Case 2: X is noninvertible and p does not divide X}}
By \newtext{Theorem~\ref{nf}}, we may take
\[
 X=\begin{pmatrix}0&1\\pa_1&b_0+pb_1\end{pmatrix},
 \qquad a_1,b_0,b_1\in K.
\]
\begin{lem}\label{lem:221-singular-first}
Every matrix in this case can be reduced to one of
\[
 \begin{pmatrix}0&1&0\\pa&b_0+pb_1&1\\u&0&w\end{pmatrix},\qquad
 \begin{pmatrix}0&1&0\\pa&b_0+pb_1&0\\1&0&w\end{pmatrix},\qquad
 \begin{pmatrix}0&1&0\\pa&b_0+pb_1&0\\0&0&w\end{pmatrix}.
\]
All parameters in this display are in $K$.
\end{lem}
\begin{proof}
Start with
\[
 E=\begin{pmatrix}0&1&v_1\\pa_1&b_0+pb_1&v_2\\u_1&u_2&w\end{pmatrix}.
\]
\textbf{Case 2.1: $\bar{v}\notin K(1,b_0)^t$.}
Put $h=v_2-b_0v_1$, so $h\neq0$. Use the change of basis
\[
 \newtext{g}=\begin{pmatrix}1&0&0\\0&1&pv_1\\-u_2h&0&h\end{pmatrix}.
\]
The matching equations give
\[
 \newtext{g}\cdot E=
 \begin{pmatrix}
 0&1&0\\
 p(a_1+v_1u_1+u_2h)&b_0+p(b_1+v_1u_2)&1\\
 u_1h&0&w-v_1u_2
 \end{pmatrix}.
\]
\edittext{Since $h\neq0$, this change of basis is invertible.}

\textbf{Case 2.2: $\bar{v}\in K(1,b_0)^t$.}
Here $v_2=b_0v_1$. If $u_1\neq0$, take
\[
 \newtext{g}=\begin{pmatrix}1&0&0\\0&1&pv_1\\-u_2/u_1&0&1/u_1\end{pmatrix}.
\]
Then
\[
 \newtext{g}\cdot E=
 \begin{pmatrix}
 0&1&0\\p(a_1+v_1u_1)&b_0+p(b_1+v_1u_2)&0\\1&0&w-v_1u_2
 \end{pmatrix}.
\]
If $u_1=0$, use
\[
 \newtext{g}=\begin{pmatrix}1&0&0\\0&1&pv_1\\-u_2&0&1\end{pmatrix},
 \qquad
 \newtext{g}\cdot E=
 \begin{pmatrix}0&1&0\\pa_1&b_0+p(b_1+v_1u_2)&0\\0&0&w-v_1u_2\end{pmatrix}.
\]
\proofreadchange{Both formulas follow from substitution in \eqref{me1}--\eqref{me4}.}
\end{proof}

\begin{lem}\label{lem:221-singular-second}
The first two forms in Lemma~\ref{lem:221-singular-first} can be reduced further by removing the coefficient $b_1$. \proofreadchange{They give types 10 and 11 in Theorem~\ref{thm:221}. The third form gives type 9.}
\end{lem}
\begin{proof}
For the first form use
\[
 \newtext{g}=\begin{pmatrix}1&0&0\\0&1&0\\-b_1b_0&b_1&1\end{pmatrix}.
\]
For the second form use
\[
 \newtext{g}=\begin{pmatrix}1&0&-pb_1\\-pb_1&1&0\\0&0&1\end{pmatrix}.
\]
In both cases the parameters change as
\[
 (a,b_1,b_0,w)\longmapsto(a+b_1b_0,0,b_0,w+b_1),
\]
and the normalized off-diagonal vectors are unchanged. \proofreadchange{In the third form both vectors vanish, so $b_1$ remains a parameter.}
\end{proof}

\subsection{\texorpdfstring{Case 3: $p\mid X$}{Case 3: p divides X}}
Write $X=pT$. \newtext{By Theorem~\ref{nf}, we may take $T=C(a,b)$ or $T=dI_2$.} \proofreadchange{Equation~\eqref{me1} modulo $p$ also gives $p\mid X'$, so we write $X'=pT'$.} Equations~\eqref{me1}--\eqref{me4} become
\begin{align}
 A_0T+\bar{b} \bar{u}^{\,t}&=T'A_0+\newtext{\bar{v}'}\bar{c}^{\,t},\label{mee1}\\
 A_0\bar{v}&=e\newtext{\bar{v}'},\label{mee2}\\
 e\bar{u}^{\,t}&=\newtext{\bar{u}'^{\,t}}A_0,\label{mee3}\\
 \bar{c}^{\,t} \bar{v}+ew&=\newtext{\bar{u}'^{\,t}}\bar{b}+ew',\label{mee4}
\end{align}
where $A_0$ denotes the reduction of $A$.

\begin{lem}\label{lem:221-pair}
In Case 3, the scalar $\bar{u}^{\,t}\bar{v}$ and the vanishing or nonvanishing of each of $\bar{u}^{\,t},\bar{v}$ are invariant under the action.
\end{lem}
\begin{proof}
Equations~\eqref{mee2} and \eqref{mee3} give $\bar{v}=A_0^{-1}e\newtext{\bar{v}'}$ and $\bar{u}^{\,t}=e^{-1}\newtext{\bar{u}'^{\,t}}A_0$. Hence $\bar{u}^{\,t}\bar{v}=\newtext{\bar{u}'^{\,t}}\newtext{\bar{v}'}$. The vanishing assertions follow from invertibility of $A_0$ and $e$.
\end{proof}

\begin{lem}[Reduction when $T$ is cyclic]\label{lem:221-cyclic-start}
Suppose that $T=C(a,b)$. The matrix $E$ can be reduced to one with
\[
 (T',\newtext{\bar{v}'},\newtext{\bar{u}'^{\,t}})=(C(a',b'),0,0),\quad
 (C(a',b'),e_1,ce_1^t),\quad
 (C(a',b'),e_1,e_2^t),\quad
 (C(a',b'),0,e_1^t),
\]
\newtext{or $(T',\newtext{\bar{v}'},\newtext{\bar{u}'^{\,t}})=(d'I_2,0,e_1^t)$. The last possibility must also be included when $\bar{v}=0$ and $\bar{u}^{\,t}\neq0$.}
\end{lem}
\begin{proof}
First suppose that $c=u_1v_1+u_2v_2\neq0$. Use the change of basis
\[
 \newtext{g}=\begin{pmatrix}
 u_1/c&u_2/c&0\\
 -v_2&v_1&0\\
 \newtext{v_2+au_2/c}&u_1/c+bu_2/c-v_1&1
 \end{pmatrix}.
\]
\proofreadchange{The upper-left block of $g$ has determinant $1$. The transformed extension has} $\newtext{\bar{v}'}=e_1$, $\newtext{\bar{u}'^{\,t}}=ce_1^t$, and $T'=C(a',b')$, where
\[
 \begin{pmatrix}u_1/c&-v_2\\u_2/c&v_1\end{pmatrix}
 \binom{a'}{b'}=\binom{av_1}{bv_1-v_2}.
\]
The coefficient matrix has determinant $1$, so $a',b'$ are uniquely determined at this step. \newtext{The last equation gives}
\[
 \newtext{w'=w+v_1(v_2+au_2/c)+v_2(u_1/c+bu_2/c-v_1).}
\]

Suppose next that $\bar{u}^{\,t}\bar{v}=0$ and both vectors are nonzero. Choose a row $\ell$ with $\ell \bar{v}=1$, and put
\[
 A_0=\begin{pmatrix}\ell\\\bar{u}^{\,t}\end{pmatrix}.
\]
\newtext{The rows are independent because $\bar{u}^{\,t}\bar{v}=0$ and $\bar{u}^{\,t}\neq0$. A block-diagonal change of basis therefore gives $\bar{v}=e_1$ and $\bar{u}^{\,t}=e_2^t$.}
Write the resulting upper-left coefficient matrix as
\[
 T=\begin{pmatrix}r&s\\\newtext{\eta}&d\end{pmatrix}.
\]
\newtext{Now take $A_0=I_2$, $e=1$, $\bar{b}=0$ and $\bar{c}^{\,t}=(r,s-1)$. Equations~\eqref{mee1}--\eqref{mee4} give $T'=C(\newtext{\eta},d)$ and $w'=w+r$. The vectors remain fixed, giving the third form.}

If $\bar{v}\neq0$ and $\bar{u}^{\,t}=0$, first choose $A_0$ with $A_0\bar{v}=e_1$. \newtext{The same operation $\bar{c}^{\,t}=(r,s-1)$ changes the first row of $T$ to $(0,1)$ and gives the second form with $c=0$.}

If $\bar{v}=0$ and $\bar{u}^{\,t}\neq0$, first normalize $\bar{u}^{\,t}=e_1^t$. \newtext{Write the resulting $T$ as above. If $s\neq0$, choose $\bar{b}=(-r,-\newtext{\eta})^t$, $\bar{c}^{\,t}=0$, $A_0=I_2$, $e=1$. This gives $T'=\left(\begin{smallmatrix}0&s\\0&d\end{smallmatrix}\right)$. The subsequent block change $A_0=\operatorname{diag}(1,s)$ gives $C(0,d)$ and leaves $\bar{u}^{\,t}=e_1^t$. If $s=0$, choose instead $\bar{b}=(d-r,-\newtext{\eta})^t$. This gives $T'=dI_2$ and $w'=w+r-d$. It is the additional scalar alternative in the statement.}

Finally, if $\bar{u}^{\,t}=\bar{v}=0$, no off-diagonal reduction is needed.
\end{proof}

\begin{lem}[Reduction when $T$ is scalar]\label{lem:221-scalar-start}
Suppose $T=dI_2$. The off-diagonal pair can be reduced to
\[
 (\bar{v},\bar{u}^{\,t})=(0,0),\quad(e_1,ce_1^t),\quad(e_1,e_2^t),\quad(0,e_1^t),
\]
while leaving $T$ scalar. The parameter $c$ may be zero.
\end{lem}
\begin{proof}
A block-diagonal change of basis leaves $dI_2$ unchanged. If $\bar{u}^{\,t}\bar{v}=c\neq0$, choose a basis with $\bar{v}=e_1$ and $\bar{u}^{\,t}=ce_1^t$. If both vectors are nonzero and $\bar{u}^{\,t}\bar{v}=0$, use the matrix with rows $\ell,\bar{u}^{\,t}$ from the preceding proof. If exactly one vector is nonzero, normalize it to $e_1$ or $e_1^t$. The remaining case is $\bar{u}^{\,t}=\bar{v}=0$.
\end{proof}

\begin{lem}[Further reductions in Case 3]\label{lem:221-identifications}
The preceding forms reduce to types 1--8 in Theorem~\ref{thm:221}.
\end{lem}
\begin{proof}
For this proof only, write
\[
 [T;\bar{v},\bar{u}^{\,t};w]=\begin{pmatrix}pT&\bar{v}\\\bar{u}^{\,t}&w\end{pmatrix}.
\]
\newtext{We use the following six reductions. Each displayed matrix $\newtext{g}$ is invertible under the stated condition, and each identity follows by substitution in \eqref{mee1}--\eqref{mee4}.}
\begin{enumerate}[label=(\arabic*)]
\item $[C(a,b);0,e_1^t;w]$ is taken to $[C(0,0);0,e_1^t;\newtext{w+b}]$ by
\[
 \newtext{g}=\begin{pmatrix}1&0&-pb\\-b&1&-pa\\0&0&1\end{pmatrix}.
\]
\item $[dI_2;e_1,e_2^t;w]$ is taken to $[C(0,0);e_1,e_2^t;w+2d]$ by
\[
 \newtext{g}=\begin{pmatrix}1&0&p\\0&1&-pd\\d&0&1\end{pmatrix}.
\]
\item $[C(a,b);e_1,e_2^t;w]$ is taken to $[C(a,0);e_1,e_2^t;w+b]$ by
\[
 \newtext{g}=\begin{pmatrix}1&0&\newtext{0}\\0&1&-pb\\0&0&1\end{pmatrix}.
\]
\item $[dI_2;e_1,ce_1^t;w]$ is taken to $[C(0,d);e_1,ce_1^t;w+d]$ by
\[
 \newtext{g}=\begin{pmatrix}1&0&0\\0&1&0\\d&-1&1\end{pmatrix}.
\]
This identity also holds when $c=0$.
\item If $c\neq0$, $[C(a,b);e_1,ce_1^t;w]$ is taken to $[C(0,b);e_1,ce_1^t;w]$ by
\[
 \newtext{g}=\begin{pmatrix}1&0&0\\0&1&-pa/c\\0&0&1\end{pmatrix}.
\]
\item If $a\neq0$, $[C(a,b);e_1,0;w]$ is taken to $[C(1,0);e_1,0;w+b]$ by
\[
 \newtext{g}=\begin{pmatrix}a&b&0\\0&1&0\\ab&a+b^2-1&a\end{pmatrix}.
\]
If $a=0$, the starting form already has type 4.
\end{enumerate}
\newtext{For completeness, the scalar form with $\bar{v}=0$, $\bar{u}^{\,t}=e_1^t$ gives type 1; item (1) gives type 2. Items (2) and (3) give type 3. Items (4)--(6) give types 4--6. The two block-diagonal forms give types 7 and 8. Thus every form from the two preceding lemmas has been accounted for.}
\end{proof}

\Needspace{24\baselineskip}
Collecting Cases 1--3 gives the following list.

\begin{thm}\label{thm:221}
\edittext{The matrices $E=\left(\begin{smallmatrix}X&\bar v\\\bar u^{\,t}&w\end{smallmatrix}\right)$ with the data below form a complete set of distinct orbit representatives. All scalar parameters, including $w$, vary independently over $K$, subject to the stated restrictions.}
\begin{center}\small
\begin{tabular}{r c c c l}
\toprule
Type & $X$ & $\bar{v}$ & $\bar{u}^{\,t}$ & Restriction\\\midrule
1 & ${p} dI_2$ & $0$ & $e_1^t$ &\\
2 & ${p} C(0,0)$ & $0$ & $e_1^t$ &\\
3 & ${p} C(a,0)$ & $e_1$ & $e_2^t$ &\\
4 & ${p} C(0,b)$ & $e_1$ & $0$ &\\
5 & ${p} C(1,0)$ & $e_1$ & $0$ &\\
6 & ${p} C(0,b)$ & $e_1$ & $c e_1^t$ & $c\neq0$\\
7 & ${p} C(a,b)$ & $0$ & $0$ &\\
8 & ${p} dI_2$ & $0$ & $0$ &\\
9 & $C({p} a,b_0+{p} b_1)$ & $0$ & $0$ &\\
10 & $C({p} a,b)$ & $e_2$ & $u e_1^t$ &\\
11 & $C({p} a,b)$ & $0$ & $e_1^t$ &\\
12 & $(d_0+{p} d_1)I_2$ & $0$ & $0$ & $d_0\neq0$\\
13 & $C(a_0+{p} a,b_0+{p} b)$ & $0$ & $0$ & $a_0\neq0$\\
14 & $dI_2+{p} C(a,b)$ & $0$ & $0$ & $d\neq0$\\\bottomrule
\end{tabular}
\end{center}
\end{thm}
\begin{proof}
\begin{blueblock}
For the coefficient comparisons and stabilizer calculations, write
\[
 X=X_0+pT,\qquad A=A_0+pL.
\]
The matching equations are equivalently
\begin{align}
 A_0X_0&=X_0'A_0,\label{me1r}\\
 A_0T+LX_0+\bar{b} \bar{u}^{\,t}&=T'A_0+X_0'L+\newtext{\bar{v}'}\bar{c}^{\,t},\label{me2r}\\
 A_0\bar{v}&=X_0'\bar{b}+e\newtext{\bar{v}'},\label{me3r}\\
 \bar{c}^{\,t} X_0+e\bar{u}^{\,t}&=\newtext{\bar{u}'^{\,t}}A_0,\label{me4r}\\
 \bar{c}^{\,t} \bar{v}+ew&=\newtext{\bar{u}'^{\,t}}\bar{b}+ew'.\label{me5r}
\end{align}
\end{blueblock}

\textbf{Comparison of representatives.}
The preceding reductions show that every orbit meets the list. Suppose that $gE=E'g^\#$. Reducing \eqref{me1} modulo $p$ gives $A_0X_0=X_0'A_0$. Hence the three cases in the reduction proof cannot meet the same orbit: $X_0$ is respectively invertible, nonzero singular, or zero. We compare the representatives within each case.

\textbf{Case 1.}
For types 12--14, both off-diagonal vectors vanish. Since $X$ and $X'$ are invertible, \eqref{me2} and \eqref{me3} force $\bar b=\bar c=0$. Thus \eqref{me1} gives ordinary similarity of $X$ and $X'$, and \eqref{me4} gives $w=w'$. Theorem~\ref{nf} determines the type and all parameters.

\textbf{Case 2.}
These are types 9--11. A comparison of types 9 and 10 in \eqref{me2} would give $e=0$. Comparing types 11 and 9 in \eqref{me3} gives the same contradiction. Types 11 and 10 are likewise separated by \eqref{me2}. Within type 9, the first and last matching equations give ordinary similarity of the companion blocks and $w=w'$, so all parameters agree.

For type 10, write $A=(a_{ij})$. The first matching equation modulo $p$ shows that $a_{21}\equiv0\pmod p$ and $b=b'$. Invertibility gives $a_{11},a_{22}\not\equiv0\pmod p$. The off-diagonal equations then give $a_{11}\equiv e\pmod p$ and $u=u'$.
For the remaining coefficient comparison, put
\[
 A=A_0+pL,\qquad
 A_0=\begin{pmatrix}e&h\\0&e+bh\end{pmatrix},\qquad e(e+bh)\neq0,
 \quad \bar b=(b_1,b_2)^t,\quad\bar c=(c_1,c_2)^t.
\]
Then $b_2=h$. The $(1,1)$ and $(2,2)$ entries of \eqref{me2r}, and \eqref{me4r}, give
\[
 L_{21}=ha+ub_1,\qquad c_2=h(a-a')+ub_1,
 \qquad c_1=uh-bc_2.
\]
Substitution in the $(2,1)$ entry of \eqref{me2r} yields
\[
 (e+bh)(a-a')=0.
\]
Since $e+bh\neq0$, we obtain $a=a'$. Consequently $c_2=ub_1$, and the last matching equation gives $w=w'$.

\begin{blueblock}
Within type 11, the reduction gives $b=b'$, and $A_0$ again has the displayed form, while \eqref{me3r} gives $b_2=0$. The $(2,2)$ and $(1,1)$ entries of \eqref{me2r} give
\[
 L_{21}=ha',\qquad b_1=h(a'-a).
\]
Its $(2,1)$ entry gives $(e+bh)(a-a')=0$. Thus $a=a'$, $b_1=0$, and \eqref{me5r} gives $w=w'$.
\end{blueblock}

\textbf{Case 3.}
Lemma~\ref{lem:221-pair} separates the cases according to the vanishing of $\bar v,\bar u$ and the value of $\bar u^t\bar v$. If both vectors vanish, ordinary similarity of $T$ separates types 7 and 8 and determines their parameters; the last matching equation gives $w=w'$.

If $\bar v=0$ and $\bar u^t=e_1^t$, the possibilities are types 1 and 2. The off-diagonal equations give $A_0=\left(\begin{smallmatrix}e&0\\h&d_*\end{smallmatrix}\right)$ with $ed_*\neq0$. Comparing a scalar $T=dI_2$ with $T'=C(0,0)$ in \eqref{mee1}, its $(1,2)$ entry would give $d_*=0$. Within type 1, the $(2,2)$ entry gives $d=d'$, and the $(1,1)$ entry then gives $b_1=0$. \newtext{Within type 2, the $(1,2)$, $(2,2)$ and $(1,1)$ entries give $d_*=e$, $h=0$ and $b_1=0$.} In both types, \eqref{mee4} gives $w=w'$.

\begin{blueblock}
For type 3, $A_0=\left(\begin{smallmatrix}e&h\\0&e\end{smallmatrix}\right)$. The $(2,1)$ entry of \eqref{mee1} gives $ea=ea'$, and its diagonal entries give $c_1=ha$ and $b_2=ha'$. Thus \eqref{mee4} gives $w=w'$.
\end{blueblock}

\begin{blueblock}
For types 4 and 5, write $T=C(a,b)$, with $(a,b)=(0,b)$ or $(1,0)$, respectively. Since $\bar v=e_1$ and $\bar u=0$, we have $A_0=\left(\begin{smallmatrix}e&h\\0&d_*\end{smallmatrix}\right)$, where $ed_*\neq0$. The $(2,1)$ entry of \eqref{mee1} gives $d_*a=a'e$, separating the two types. Within type 4, the $(2,2)$ and $(1,1)$ entries give $b=b'$ and $c_1=0$. Within type 5, they give $h=0$ and $c_1=0$. Hence \eqref{mee4} gives $w=w'$ in both cases.

In type 6, the invariant $\bar u^t\bar v$ gives $c=c'\neq0$, and $A_0=\operatorname{diag}(e,d_*)$. The $(2,2)$ entry gives $b=b'$, while the $(1,1)$ entry gives $c_1=cb_1$. Equation~\eqref{mee4} then gives $w=w'$.
\end{blueblock}
Thus the families are pairwise disjoint, and no family has repeated parameter choices.

\end{proof}

\subsection{Explicit stabilizers}
\begin{blueblock}
\edittext{Set $E'=E$ in \eqref{me1r}--\eqref{me5r}. In types 1--8, $L\in M_2(K)$ is arbitrary. The remaining free parameters lie in $K$, with $e\neq0$ and the further restrictions shown in the table.}
\end{blueblock}
\begin{center}\small
\begin{tabular}{r c c c l}
\toprule
Type & $A_0$ & $\bar{b}$ & $\bar{c}^{\,t}$ & Condition\\\midrule
1 & $\begin{pmatrix}e&0\\h&d_*\end{pmatrix}$ & $0$ & $(r,s)$ & $d_*\neq0$\\[8pt]
2 & $eI_2$ & $0$ & $(r,s)$ &\\[4pt]
3 & $\begin{pmatrix}e&h\\0&e\end{pmatrix}$ & $(r,ah)^t$ & $(ah,r)$ &\\[8pt]
4 & $\begin{pmatrix}e&h\\0&d\end{pmatrix}$ & $(r,s)^t$ & $(0,e+bh-d)$ & $d\neq0$\\[8pt]
5 & $eI_2$ & $(r,s)^t$ & $0$ &\\[4pt]
6 & $\operatorname{diag}(e,d)$ & $(r,0)^t$ & $(cr,e-d)$ & $d\neq0$\\[4pt]
7 & $rI_2+sC(a,b)$ & arbitrary & arbitrary & $\det A_0\neq0$\\[4pt]
8 & arbitrary & arbitrary & arbitrary & $A_0\in\mathrm{GL}_2(K)$\\\bottomrule
\end{tabular}
\end{center}

For type 9 the full description is
\[
 A=rI_2+sX,\quad r,s\in S_2,\quad\det A\in S_2^\times,
 \qquad\bar{b}=(u,0)^t,\quad\bar{c}^{\,t}=(\newtext{-b_0v},v),\quad u,v\in K,
\]
with $e\in K^\times$ arbitrary. For types 10 and 11 put
\[
 A_0=\begin{pmatrix}e&h\\0&e+bh\end{pmatrix},\qquad e(e+bh)\neq0.
\]
Their remaining entries are, respectively,
\begin{align*}
\text{Type 10:}\quad
 L&=\newtext{\begin{pmatrix}\ell-bj&j\\ur+ah&\ell\end{pmatrix}},&
 \bar{b}&=(r,h)^t,&\bar{c}^{\,t}&=(u(h-br),ur),\\
\text{Type 11:}\quad
 L&=\newtext{\begin{pmatrix}\ell-bj&j\\ah&\ell\end{pmatrix}},&
 \bar{b}&=0,&\bar{c}^{\,t}&=(h-br,r),
\end{align*}
where $h,r,j,\ell\in K$. Finally, in types 12--14, $\bar{b}=\bar{c}^{\,t}=0$, $e\in K^\times$ is arbitrary, and
\[
\begin{array}{c|l}
12&A\in\mathrm{GL}_2(S_2)\text{ arbitrary},\\
13&A=rI_2+sX,\quad r,s\in S_2,\quad\det A\in S_2^\times,\\
14&A=A_0+{p} L,\quad A_0=rI_2+sC(a,b)\in\mathrm{GL}_2(K),\quad L\in M_2(K).
\end{array}
\]
\begin{blueblock}
\edittext{For type 10, the comparison proof already gives the displayed $\bar b$ and $\bar c$. The remaining entries of \eqref{me2r} give $L_{21}=ur+ah$ and $L_{11}=L_{22}-bL_{12}$. This yields the stated form of $L$. Substitution in \eqref{me1r}--\eqref{me5r} gives the other types in the same way.}

\end{blueblock}
\subsection{Counting and the class equation}
\begin{blueblock}
For $|K|=q$, let $\theta\in\{\mathrm I,\mathrm R,\mathrm S\}$ denote the irreducible, repeated-root and split-distinct-root types of a monic quadratic. Write $\newtext{\nu_\theta}$ for the number of all such quadratics, $\newtext{\nu_\theta^\times}$ for the number with nonzero constant term, and $\newtext{h_\theta}$ for the centralizer order of their companion matrices in $\mathrm{GL}_2(K)$. Then
{\color{black}
\begin{equation}\label{eq:quadratic-data}
\begin{array}{c|c|c|c}
\theta&\newtext{\nu_\theta}&\newtext{\nu_\theta^\times}&\newtext{h_\theta}\\\hline
\mathrm I&q(q-1)/2&q(q-1)/2&q^2-1\\
\mathrm R&q&q-1&q(q-1)\\
\mathrm S&q(q-1)/2&\newtext{(q-1)(q-2)/2}&(q-1)^2.
\end{array}
\end{equation}
}
The centralizer is the group of units in $K[T]/(f)$\proofreadchange{, where $f$ is the characteristic polynomial. Its factorization gives the order in each case. To count the split quadratics, choose an unordered pair of distinct roots or a repeated root. The remaining monic quadratics are irreducible.} Requiring nonzero constant term excludes zero roots. Thus these formulas hold also in characteristic two.

\end{blueblock}
\edittext{Put $G=|G_\lambda|=q^9(q-1)^2(q^2-1)$. The columns have the meanings fixed in Section~\ref{subsec:extension-action}. Each row indexed by $\theta$ contributes once for each of the three quadratic types in \eqref{eq:quadratic-data}.}
\afterpage{\begin{landscapetablepages}
\begin{center}\normalfont\small
Orbit counts and the class equation for $\lambda=(2,2,1)$

\medskip
\renewcommand{\arraystretch}{1.5}
\begin{tabular}{r l c c c c}
\toprule
Type & Condition & $N$ & $H$ & $|\mathcal O|=G/H$ & $N|\mathcal O|$\\\midrule
1 & & $q^2$ & $q^7(q-1)^2$ & $q^2(q^2-1)$ & $q^4(q^2-1)$\\
2 & & $q$ & $q^6(q-1)$ & $q^3(q-1)(q^2-1)$ & $q^4(q-1)(q^2-1)$\\
3 & & $q^2$ & $q^6(q-1)$ & $q^3(q-1)(q^2-1)$ & $q^5(q-1)(q^2-1)$\\
4 & & $q^2$ & $q^7(q-1)^2$ & $q^2(q^2-1)$ & $q^4(q^2-1)$\\
5 & & $q$ & $q^6(q-1)$ & $q^3(q-1)(q^2-1)$ & $q^4(q-1)(q^2-1)$\\
6 & $c\neq0$ & $q^2(q-1)$ & $q^5(q-1)^2$ & $q^4(q^2-1)$ & $q^6(q-1)(q^2-1)$\\
7 & $T^2-bT-a$ of type $\theta$ & $q\newtext{\nu_\theta}$ & $q^8(q-1)\newtext{h_\theta}$ & $\dfrac{q(q-1)(q^2-1)}{\newtext{h_\theta}}$ & $\dfrac{q^2(q-1)(q^2-1)\newtext{\nu_\theta}}{\newtext{h_\theta}}$\\[6pt]
8 & & $q^2$ & $G$ & $1$ & $q^2$\\
9 & $b_0=0$ & $q^3$ & $q^5(q-1)^2$ & $q^4(q^2-1)$ & $q^7(q^2-1)$\\
9 & $b_0\neq0$ & $q^3(q-1)$ & $q^4(q-1)^3$ & $q^5(q+1)$ & $q^8(q^2-1)$\\
10 & $b=0$ & $q^3$ & $q^4(q-1)$ & $q^5(q-1)(q^2-1)$ & $q^8(q-1)(q^2-1)$\\
10 & $b\neq0$ & $q^3(q-1)$ & $q^3(q-1)^2$ & $q^6(q^2-1)$ & $q^9(q-1)(q^2-1)$\\
11 & $b=0$ & $q^2$ & $q^4(q-1)$ & $q^5(q-1)(q^2-1)$ & $q^7(q-1)(q^2-1)$\\
11 & $b\neq0$ & $q^2(q-1)$ & $q^3(q-1)^2$ & $q^6(q^2-1)$ & $q^8(q-1)(q^2-1)$\\
12 & $d_0\neq0$ & $q^2(q-1)$ & $q^5(q-1)^3(q+1)$ & $q^4$ & $q^6(q-1)$\\
13 & $T^2-b_0T-a_0$ of type $\theta$ & $q^3\newtext{\nu_\theta^\times}$ & $q^2(q-1)\newtext{h_\theta}$ & $\dfrac{q^7(q-1)(q^2-1)}{\newtext{h_\theta}}$ & $\dfrac{q^{10}(q-1)(q^2-1)\newtext{\nu_\theta^\times}}{\newtext{h_\theta}}$\\[6pt]
14 & $T^2-bT-a$ of type $\theta$ & $q(q-1)\newtext{\nu_\theta}$ & $q^4(q-1)\newtext{h_\theta}$ & $\dfrac{q^5(q-1)(q^2-1)}{\newtext{h_\theta}}$ & $\dfrac{q^6(q-1)^2(q^2-1)\newtext{\nu_\theta}}{\newtext{h_\theta}}$\\[6pt]\bottomrule
\end{tabular}
\end{center}
\end{landscapetablepages}}\par
\begin{blueblock}
The stabilizer orders follow directly from the free parameters above. For types 1--8, the arbitrary $L$ contributes $q^4$. For types 10 and 11, the two conditions $e\neq0$ and $e+bh\neq0$ give $q(q-1)$ choices for $(e,h)$ if $b=0$, and $(q-1)^2$ if $b\neq0$. The other three free parameters contribute $q^3$. The centralizer descriptions for types 9 and 12--14 give the remaining rows.

\end{blueblock}
\begin{blueblock}
\Needspace{14\baselineskip}
\proofreadchange{Summing the last column over all rows with zero, nonzero singular and invertible reduction gives, respectively,}
\begin{align*}
 \proofreadchange{\sum_{\text{types }1\text{--}8}NG/H}&=q^9,\\*
 \proofreadchange{\sum_{\text{types }9\text{--}11}NG/H}&=q^9(q^3+q^2-q-1),\\*
 \proofreadchange{\sum_{\text{types }12\text{--}14}NG/H}&=q^{10}(q-1)^2(q+1).
\end{align*}
Their sum is $q^{13}=|E_{(2,2,1)}|$. This also agrees with counting the possible reductions $X_0$: there are $1$, $q^3+q^2-q-1$ and $|\mathrm{GL}_2(K)|$ choices in the three cases, each with $q^9$ choices of the remaining data. Summing the column $N$ gives
\begin{equation}\label{eq:c221}
 c_{(2,2,1)}(q)=q^5+2q^4+3q^3+2q^2+2q.
\end{equation}

\end{blueblock}
\clearpage
\section{The family \texorpdfstring{$(m+1,1,1,1)$}{(m+1,1,1,1)}}\label{sec:hook}
\proofreadchange{Let $K$ be any field and $m\geq1$. In this section, $p$ is an indeterminate and $S_j=K[p]/(p^j)$.} For this partition,
\[
 M_\lambda=S_{m+1}\oplus K^3,
\]
and the extension matrices and automorphisms have the forms
\[
 E=\begin{pmatrix}x&\newtext{\bar{y}}^t\\\newtext{\bar{z}}&W\end{pmatrix},\qquad
 g=\begin{pmatrix}a&p^m\newtext{\bar{b}^{\,t}}\\\newtext{\bar{c}}&D\end{pmatrix},
\]
where $x\in S_{m+1}$, $a\in S_{m+1}^\times$, $\newtext{\bar y,\bar z,\bar b,\bar c\in K^3}$, $W\in M_3(K)$ and $D\in\mathrm{GL}_3(K)$. \newtext{Thus $\bar y^t$ and $\bar b^t$ are rows, while $\bar z$ and $\bar c$ are columns. The unbarred letter $z$ below denotes the scalar parameter in the normal forms.} \proofreadchange{Let $a_0,x_0\in K$ satisfy $a_0\equiv a\pmod{p}$ and $x_0\equiv x\pmod{p}$.}

\begin{lem}\label{lem:hook-pair}
If $E$ and $E'$ are in the same orbit and $p\mid x$, then $p\mid x'$, $\newtext{\bar{y}}^t\newtext{\bar{z}}=\newtext{\bar{y}'}^t\newtext{\bar{z}'}$, and
\[
 \newtext{\bar{y}}=0\Longleftrightarrow \newtext{\bar{y}'}=0,\qquad \newtext{\bar{z}}=0\Longleftrightarrow \newtext{\bar{z}'}=0.
\]
\end{lem}
\begin{proof}
The upper-left matching equation modulo $p$ gives $a_0x_0=a_0x'_0$, so $x_0=x'_0$. When these residues are zero, the off-diagonal equations become $a_0\newtext{\bar{y}}^t=\newtext{\bar{y}'}^tD$ and $D\newtext{\bar{z}}=\newtext{\bar{z}'}a_0$. Multiplying gives $\newtext{\bar{y}}^t\newtext{\bar{z}}=\newtext{\bar{y}'}^t\newtext{\bar{z}'}$. Since $a_0$ and $D$ are invertible, the two vanishing assertions follow.
\end{proof}

\edittext{Let $e_1,e_2,e_3$ be the standard basis of $K^3$, and let $J_r$ be the nilpotent Jordan block with ones immediately above the diagonal.} As before, $C(w_1,w_2)=\left(\begin{smallmatrix}0&1\\w_1&w_2\end{smallmatrix}\right)$. \proofreadchange{Fix a complete set $\mathcal R_3(K)$ of rational canonical representatives in $M_3(K)$.}
\Needspace{30\baselineskip}
\begin{thm}\label{thm:hook}
The following table gives a complete set of distinct representatives $E=\left(\begin{smallmatrix}x&\newtext{\bar{y}}^t\\\newtext{\bar{z}}&W\end{smallmatrix}\right)$. In types 1--12, $x\in pS_{m+1}$. All other scalar parameters are in $K$; $\newtext{z}$ may be zero in types 1 and 2, and $w_1\neq0$ in type 6. In type 13, $x\in S_{m+1}$ is arbitrary.
\begin{center}\small
\begin{tabular}{r c c c}
\toprule
Type & $\newtext{\bar{y}}^t$ & $\newtext{\bar{z}}$ & $W$\\\midrule
1 & $e_1^t$ & $\newtext{z} e_1$ & $0\oplus wI_2$\\
2 & $e_1^t$ & $\newtext{z} e_1$ & $0\oplus C(w_1,w_2)$\\[3pt]
3 & $e_1^t$ & $e_2$ & $\begin{pmatrix}0&0&1\\0&0&0\\0&w&0\end{pmatrix}$\\[12pt]
4 & $e_1^t$ & $e_2$ & $\begin{pmatrix}0&0&0\\0&0&0\\0&1&0\end{pmatrix}$\\[12pt]
5 & $e_1^t$ & $e_2$ & $\operatorname{diag}(0,0,w)$\\[3pt]
6 & $e_1^t$ & $e_2$ & $\begin{pmatrix}0&w_1&0\\0&0&0\\0&0&w_2\end{pmatrix}$\\[12pt]
7 & $e_1^t$ & $0$ & $J_3$\\
8 & $e_1^t$ & $0$ & $J_2\oplus w$\\
9 & $0$ & $e_1$ & $0\oplus wI_2$\\
10 & $0$ & $e_1$ & $0\oplus C(w_1,w_2)$\\
11 & $0$ & $e_1$ & $J_3^t$\\
12 & $0$ & $e_1$ & $J_2^t\oplus w$\\
13 & $0$ & $0$ & $W\in\mathcal R_3(K)$\\\bottomrule
\end{tabular}
\end{center}
\end{thm}

\begin{proof}
We follow the cases determined by $x,\newtext{\bar{y}},\newtext{\bar{z}}$.

\textbf{Case 1: $p\nmid x$.}
Proposition~\ref{1stt} reduces $E$ to a block-diagonal matrix. A further change by $\operatorname{diag}(1,D)$ puts the lower block into rational canonical form. This gives type 13 with $x$ a unit.

\textbf{Case 2: $p\mid x$.}
For two matrices in this case, $gE=E'g^\#$ gives
\begin{align}
 a(x-x')&\mathrel{\proofreadchange{\equiv}}p^m(\newtext{\bar{y}'}^t\newtext{\bar{c}}-\newtext{\bar{b}^{\,t}} \newtext{\bar{z}})\pmod{p^{m+1}},\label{eq:hook1}\\
 a_0\newtext{\bar{y}}^t&=\newtext{\bar{y}'}^tD,\qquad D\newtext{\bar{z}}=\newtext{\bar{z}'}a_0,\label{eq:hook2}\\
 \newtext{\bar{c}} \newtext{\bar{y}}^t+DW&=W'D+\newtext{\bar{z}'}\newtext{\bar{b}^{\,t}}.\label{eq:hook3}
\end{align}
Thus
\begin{equation}\label{eq:hook-change-W}
 W'=DWD^{-1}+(\newtext{\bar{c}} \newtext{\bar{y}}^t-\newtext{\bar{z}'}\newtext{\bar{b}^{\,t}})D^{-1}.
\end{equation}
\newtext{Once \proofreadchange{$a$}, $\newtext{\bar{b}^{\,t}},\newtext{\bar{c}},D$ have been chosen, \eqref{eq:hook1} determines $x'$. It changes only the coefficient of $p^m$, so $x'$ is still in $pS_{m+1}$. We may therefore carry out the following reductions of $W$ without imposing any additional restriction on $x$.}

\textbf{Case 2.1: $\newtext{\bar{y}}\neq0$ and $\newtext{\bar{z}}\neq0$.}
First choose $D$ so that $\newtext{\bar{y}'}^t=e_1^t$. The matrices preserving this row have the form
\[
 D=\begin{pmatrix}a_0&0&0\\u_1&d_1&d_2\\u_2&d_3&d_4\end{pmatrix},
 \qquad \begin{pmatrix}d_1&d_2\\d_3&d_4\end{pmatrix}\in\mathrm{GL}_2(K).
\]
Write $\newtext{\bar{z}=\binom{z_1}{\proofreadchange{\bar z_2}}}$, where $\proofreadchange{\bar z_2}\in K^2$. If $z_1\neq0$, choose $u_1,u_2$ to make the last two entries of $\newtext{\bar{z}'}$ zero. If $z_1=0$, \proofreadchange{take $a_0=1$ and} choose the lower two by two block to send $\proofreadchange{\bar z_2}$ to $(1,0)^t$. Hence the two possibilities are
\[
 (\newtext{\bar{y}}^t,\newtext{\bar{z}})=(e_1^t,\newtext{z} e_1),\quad\newtext{z}\neq0,
 \qquad\text{or}\qquad (\newtext{\bar{y}}^t,\newtext{\bar{z}})=(e_1^t,e_2).
\]

Suppose first that $\newtext{\bar{z}}=\newtext{z} e_1$, with $\newtext{z}\neq0$. Equation~\eqref{eq:hook-change-W} allows us to make the first row and column of $W$ zero. \newtext{For example, take $D=I_3$, $a=1$, $\newtext{b_1}=0$, $\newtext{c_i}=-w_{i1}$, and $\newtext{b_j}=w_{1j}/\newtext{z}$ for $j=2,3$.} The remaining lower two by two block is reduced by ordinary similarity. It is either scalar or a companion matrix. We obtain types 1 and 2 with $\newtext{z}\neq0$.

Now suppose that $\newtext{\bar{z}}=e_2$. Choose $\newtext{\bar{c}},\newtext{\bar{b}^{\,t}}$ to make the first column and second row of $W$ zero. \newtext{One choice is $D=I_3$, $a=1$, $\newtext{c_i}=-w_{i1}$ and $\newtext{\bar{b}^{\,t}}=(0,w_{22},w_{23})$.} Write the remaining matrix as
\[
 W=\begin{pmatrix}0&\newtext{w_{12}}&\newtext{w_{13}}\\0&0&0\\0&\newtext{w_{32}}&\newtext{w_{33}}\end{pmatrix}.
\]
\begin{blueblock}
The further basis changes preserving $\newtext{\bar{y}}^t=e_1^t$ and $\newtext{\bar{z}}=e_2$ have
\[
 D=\begin{pmatrix}r&0&0\\u&r&v\\s&0&\newtext{\eta}\end{pmatrix},\qquad \proofreadchange{r=a_0},\quad r\newtext{\eta}\neq0.
\]
The first column and second row can again be made zero by $\newtext{\bar{c}},\newtext{\bar{b}^{\,t}}$. Comparing the other four entries in \eqref{eq:hook3} gives
\begin{equation}\label{eq:hook-residual}
 \newtext{w_{12}'}=\newtext{w_{12}},\quad r\newtext{w_{13}}=\newtext{w_{12}} v+\newtext{w_{13}'}\newtext{\eta},\quad
 s\newtext{w_{12}}+\newtext{\eta}\newtext{w_{32}}=\newtext{w_{32}'}r,\quad s\newtext{w_{13}}+\newtext{\eta}\newtext{w_{33}}=\newtext{w_{32}'}v+\newtext{w_{33}'}\newtext{\eta}.
\end{equation}
\edittext{We now distinguish four cases according to the entries of $W$.}
\end{blueblock}

If $\newtext{w_{12}}\neq0$, \newtext{choose $r=\newtext{\eta}=1$, $v=\newtext{w_{13}}/\newtext{w_{12}}$, $s=-\newtext{w_{32}}/\newtext{w_{12}}$. This gives $\newtext{w_{13}'}=\newtext{w_{32}'}=0$ and $\newtext{w_{33}'}=\newtext{w_{33}}-\newtext{w_{13}}\newtext{w_{32}}/\newtext{w_{12}}$.} We obtain
\[
 W'=\begin{pmatrix}0&w_1&0\\0&0&0\\0&0&w_2\end{pmatrix},\qquad w_1\neq0,
\]
and hence type 6.

If $w_{12}=0$ and $\newtext{w_{13}}\neq0$, \newtext{choose $r=1$, $\newtext{\eta}=\newtext{w_{13}}$, $v=0$ and $s=-\newtext{w_{33}}$. Then $\newtext{w_{13}'}=1$, $\newtext{w_{32}'}=\newtext{w_{13}}\newtext{w_{32}}$ and $\newtext{w_{33}'}=0$.} This gives
\[
 W'=\begin{pmatrix}0&0&1\\0&0&0\\0&w&0\end{pmatrix},
\]
and hence type 3.

If $w_{12}=w_{13}=0$ and $\newtext{w_{32}}\neq0$, \newtext{choose $r=1$, $\newtext{\eta}=1/\newtext{w_{32}}$, $s=0$ and $v=\newtext{w_{33}}/\newtext{w_{32}}$. Then $\newtext{w_{32}'}=1$ and $\newtext{w_{33}'}=0$.} This gives type 4. If all three entries vanish, the matrix is already $\operatorname{diag}(0,0,\newtext{w_{33}})$, giving type 5.

\textbf{Case 2.2: $\newtext{\bar{y}}\neq0$ and $\newtext{\bar{z}}=0$.}
Normalize $\newtext{\bar{y}}^t=e_1^t$, and choose $\newtext{\bar{c}}$ to make the first column of $W$ zero. Write
\[
 W=\begin{pmatrix}0&\ell\\0&B\end{pmatrix},\qquad
 \ell=(w_{12},w_{13}),\quad B\in M_2(K).
\]
If $\ell=0$, reduce $B$ to a scalar or companion matrix by ordinary similarity. This gives types 1 and 2 with $\newtext{z}=0$.

Suppose that $\ell\neq0$. \newtext{For $D=\left(\begin{smallmatrix}r&0\\v&H\end{smallmatrix}\right)$, equation \eqref{eq:hook-change-W}, followed by making the first column zero again, gives}
\[
 \newtext{\ell'=r\ell H^{-1},\qquad B'=HBH^{-1}+v\ell H^{-1}.}
\]
\newtext{\proofreadchange{Take $r=1$ throughout this reduction and choose $H$ to make $\ell'=(1,0)$.} Then take $H=I_2$ and choose $v$ to make the first column of $B$ zero. We are left with}
\[
 \newtext{B=\begin{pmatrix}0&b\\0&d\end{pmatrix}.}
\]
\newtext{If $b\neq0$, use $H=\left(\begin{smallmatrix}1&0\\-d&b\end{smallmatrix}\right)$ and choose $v$ to keep the first column zero. This changes $b,d$ to $1,0$. If $b=0$, the scalar $d$ remains.} The resulting lower blocks are, respectively,
\[
 W'=\begin{pmatrix}0&1&0\\0&0&1\\0&0&0\end{pmatrix},\qquad
 W'=\begin{pmatrix}0&1&0\\0&0&0\\0&0&w\end{pmatrix}.
\]
These give types 7 and 8.

\textbf{Case 2.3: $\newtext{\bar{y}}=0$ and $\newtext{\bar{z}}\neq0$.}
Normalize $\newtext{\bar{z}}=e_1$, and choose $\newtext{\bar{b}^{\,t}}$ to make the first row of $W$ zero. If $(w_{21},w_{31})^t=0$, reduce the lower two by two block to scalar or companion form. This gives types 9 and 10. If $(w_{21},w_{31})^t\neq0$, the transpose of the reductions in Case 2.2 gives
\[
 W'=\begin{pmatrix}0&0&0\\1&0&0\\0&1&0\end{pmatrix},\qquad
 W'=\begin{pmatrix}0&0&0\\1&0&0\\0&0&w\end{pmatrix}.
\]
These give types 11 and 12. \newtext{The transposed operations are permitted: transposing $gE=E'g^\#$ replaces $g$ by $((g^\#)^t)^{-1}$, which belongs to the same module automorphism group.}

\textbf{Case 2.4: $\newtext{\bar{y}}=\newtext{\bar{z}}=0$.}
Ordinary similarity puts the lower block in $\mathcal R_3(K)$, giving type 13 also when $p\mid x$.

\textbf{Comparison of representatives.}
The residue of $x$ is unchanged under the action. Thus the unit case cannot meet the case $p\mid x$. For the block-diagonal representatives of type 13, the matching equations give $x=x'$ and $DW=W'D$. Since $W,W'\in\mathcal R_3(K)$, this implies $W=W'$. Lemma~\ref{lem:hook-pair} separates the remaining cases according to the vanishing of $\bar y,\bar z$ and their product $\bar y^t\bar z$. It remains to compare representatives within these cases, using $gE=E'g^\#$.

\begin{blueblock}
We first record the invariant needed to recover $x$. When $p\mid x$, write $x=\sum_{i=1}^m x_ip^i$ with $x_i\in K$, using the fixed coefficient representatives. Equation~\eqref{eq:hook1} and the invertibility of $a$ imply that $x-x'\in p^mS_{m+1}$. Thus $x_i=x_i'$ for $1\leq i<m$: only the highest coefficient $x_m$ can change. For this partition, define
\begin{equation}\label{eq:hook-trace}
 \tau=x_m+\operatorname{tr}W\in K.
\end{equation}
\edittext{Thus $\tau$ is the coefficient of $p^m$ in $x+p^m\operatorname{tr}W$.}

Taking the coefficient of $p^m$ in \eqref{eq:hook1} gives
\[
 x_m-x_m'=a_0^{-1}\bigl(\bar y'^{\,t}\bar c-\bar b^{\,t}\bar z\bigr).
\]
On the other hand, multiplying \eqref{eq:hook3} on the right by $D^{-1}$ and taking traces gives
\begin{align*}
 \operatorname{tr}W'-\operatorname{tr}W
 &=\bar y^{\,t}D^{-1}\bar c-\bar b^{\,t}D^{-1}\bar z'\\
 &=a_0^{-1}\bigl(\bar y'^{\,t}\bar c-\bar b^{\,t}\bar z\bigr)
 =x_m-x_m',
\end{align*}
where the second equality uses \eqref{eq:hook2}. Consequently $x_m+\operatorname{tr}W=x_m'+\operatorname{tr}W'$, so $\tau$ is constant on each orbit. In the case $(2,1,1,1)$, we have $m=1$, so $x=px_1$ and $\tau=x_1+\operatorname{tr}W$.
\end{blueblock}

\begin{blueblock}
Write $r=a_0\neq0$.

\textit{The pair $(\bar y^t,\bar z)=(e_1^t,ze_1)$, $z\neq0$.}
The scalar $z$ is invariant, and \eqref{eq:hook2} forces $D=\operatorname{diag}(r,H)$ with $H\in\mathrm{GL}_2(K)$. If $W=0\oplus B$ and $W'=0\oplus B'$, the lower two by two block of \eqref{eq:hook3} is $HB=B'H$. The scalar and companion forms for $B$ therefore distinguish types 1 and 2 and determine their parameters.

\textit{The pair $(\bar y^t,\bar z)=(e_1^t,e_2)$.}
Use \eqref{eq:hook-residual}, where $r\newtext{\eta}\neq0$. The entry $w_{12}$ is invariant, so type 6 is separated from types 3--5. In type 6, $w_{13}=w_{32}=w_{13}'=w_{32}'=0$ forces $v=s=0$, and then $w_{33}=w_{33}'$. If $w_{12}=0$, the equations preserve the vanishing of each of $w_{13},w_{32}$. Type 3 has $w_{13}=1$, and the second and third equations give $\newtext{\eta}=r$ and $w_{32}=w_{32}'$. Type 4 has $w_{13}=0,w_{32}=1$ and no free parameter in $W$. Type 5 has both entries zero, and the last equation gives $w_{33}=w_{33}'$. Thus these four types are disjoint and their parameters agree whenever they are equivalent.

\textit{The pair $(\bar y^t,\bar z)=(e_1^t,0)$.}
Write
\[
 D=\begin{pmatrix}r&0\\v&H\end{pmatrix},\qquad
 W=\begin{pmatrix}0&\ell\\0&B\end{pmatrix},\qquad H\in\mathrm{GL}_2(K).
\]
Comparing the upper-right and lower-right blocks in \eqref{eq:hook3} gives
\[
 r\ell=\ell'H,\qquad v\ell+HB=B'H.
\]
Consequently $\ell=0$ if and only if $\ell'=0$. In this case $HB=B'H$, which distinguishes types 1 and 2 with $z=0$ and determines their parameters. Otherwise both normal forms have $\ell=\ell'=(1,0)$, so
\[
 H=\begin{pmatrix}r&0\\h&\newtext{\eta}\end{pmatrix},\qquad r\newtext{\eta}\neq0.
\]
For types 7 and 8, respectively, $B=\left(\begin{smallmatrix}0&\beta\\0&\delta\end{smallmatrix}\right)$ with $(\beta,\delta)=(1,0)$ or $(0,w)$. The second-column comparison gives
\[
 r\beta=\beta'\newtext{\eta},\qquad h\beta+\newtext{\eta}\delta=\delta'\newtext{\eta}.
\]
The first equation separates types 7 and 8; within type 8 the second gives $w=w'$. Type 7 has no free parameter in $W$.

\textit{The pair $(\bar y^t,\bar z)=(0,e_1)$.}
Transpose the preceding comparison. As explained in Case 2.3, transposition preserves orbit equivalence, so it separates types 9--12 and determines their parameters.

In every case above, the displayed lower blocks satisfy $W=W'$ once their parameters have been determined. Since $\tau=\tau'$ by \eqref{eq:hook-trace}, we obtain
\[
 x_m=\tau-\operatorname{tr}W
     =\tau'-\operatorname{tr}W'=x_m'.
\]
\edittext{Together with $x_i=x_i'$ for $i<m$, this gives $x=x'$. Hence every orbit contains exactly one listed representative.}
\end{blueblock}

\end{proof}


\subsection{Explicit stabilizers}
\begin{blueblock}
In types 1, 2, 9 and 10, write $B$ for the lower two by two block: it is scalar in types 1 and 9, and companion in types 2 and 10. \edittext{In types 1--12, put $r=a_0\neq0$; the $m$ higher coefficients of $a$ are arbitrary.} The remaining parameters are in $K$. The equations for a stabilizer are
\begin{equation}\label{eq:hook-stabilizer}
 r\bar{y}^t=\bar{y}^tD,\qquad D\bar{z}=r\bar{z},\qquad
 \bar{b}^{\,t} \bar{z}=\bar{y}^t\bar{c},\qquad DW-WD+\bar{c} \bar{y}^t-\bar{z}\bar{b}^{\,t}=0.
\end{equation}
\end{blueblock}
In types $1,2$ with $z\neq0$ they give
\[
 D=\operatorname{diag}(r,H),\quad H\in C_{\mathrm{GL}_2(K)}(B),\qquad
 \bar{b}^{\,t}=(u,0,0),\quad\bar{c}=(z u,0,0)^t.
\]
\edittext{For types 3--6, where $(\bar y^t,\bar z)=(e_1^t,e_2)$, the solutions are}
\begin{center}\small
\begin{tabular}{c c c c}
\toprule
Type & $D$ & $\bar{b}^{\,t}$ & $\bar{c}^t$\\\midrule
$3$ & $\begin{pmatrix}r&0&0\\u&r&v\\wv&0&r\end{pmatrix}$ & $(\newtext{\eta},wv,u)$ & $(wv,\newtext{\eta},wu)$\\[12pt]
$4$ & $\begin{pmatrix}r&0&0\\u&r&0\\v&0&r\end{pmatrix}$ & $(\newtext{\eta},0,0)$ & $(0,\newtext{\eta},u)$\\[12pt]
$5$ & $\begin{pmatrix}r&0&0\\u&r&v\\s&0&h\end{pmatrix}$ & $(\newtext{\eta},0,wv)$ & $(0,\newtext{\eta},ws)$\\[12pt]
$6$ & $\begin{pmatrix}r&0&0\\u&r&0\\0&0&h\end{pmatrix}$ & $(\newtext{\eta},w_1 u,0)$ & $(w_1 u,\newtext{\eta},0)$\\\bottomrule
\end{tabular}
\end{center}
Here $h\neq0$ in the last two rows, and $u,v,s,\newtext{\eta}$ are free whenever present.
\edittext{When exactly one of $\bar y,\bar z$ is nonzero, the solutions are:}
\begin{align*}
1,2\ (z=0):\quad &D=\begin{pmatrix}r&0\\v&H\end{pmatrix},\quad
 H\in C_{\mathrm{GL}_2(K)}(B),\quad v\in K^2,\quad
 \bar{c}=\binom{0}{Bv},\quad\bar{b}^{\,t}\text{ arbitrary};\\
9,10:\quad &D=\begin{pmatrix}r&v\\0&H\end{pmatrix},\quad
 H\in C_{\mathrm{GL}_2(K)}(B),\quad v\in\proofreadchange{K^{1\times 2}},\quad
 \bar{b}^{\,t}=(0,vB),\quad\bar{c}\text{ arbitrary};\\
7:\quad &D=rI_3,\quad\bar{c}=0,\quad\bar{b}^{\,t}\text{ arbitrary};\\
11:\quad &D=rI_3,\quad\bar{b}^{\,t}=0,\quad\bar{c}\text{ arbitrary};\\
8:\quad &D=\begin{pmatrix}r&0&0\\0&r&0\\uw&u&h\end{pmatrix},\quad
 \bar{c}=(0,0,uw^2)^t,\quad\bar{b}^{\,t}\text{ arbitrary},\quad h\neq0;\\*
12:\quad &D=\begin{pmatrix}r&0&uw\\0&r&u\\0&0&h\end{pmatrix},\quad
 \bar{b}^{\,t}=(0,0,uw^2),\quad\bar{c}\text{ arbitrary},\quad h\neq0.
\end{align*}
\begin{blueblock}
\edittext{For types 3--6, the first two equations of \eqref{eq:hook-stabilizer} give $D=\left(\begin{smallmatrix}r&0&0\\u&r&v\\s&0&h\end{smallmatrix}\right)$. Substituting each $W$ in the remaining equations gives the table. For types 1 and 2 with $z=0$, the same calculation gives $HB=BH$ and $\bar c=\binom{0}{Bv}$. The other cases follow by substitution in \eqref{eq:hook-stabilizer}; no division by $w$ is required.}

\end{blueblock}
\proofreadchange{For a block-diagonal representative of type 13, the stabilizer is}
\[
\proofreadchange{\Stab_{G_\lambda}(E)}=
\begin{cases}
\{\operatorname{diag}(a,D):a\in \reviewtext{S_{m+1}^\times},\ D\in C_{\mathrm{GL}_3(K)}(W)\},&x\in \reviewtext{S_{m+1}^\times},\\[2pt]
\left\{\begin{pmatrix}a&{p}^m\bar{b}^{\,t}\\\bar{c}&D\end{pmatrix}:a\in \reviewtext{S_{m+1}^\times},\ \bar{b}^{\,t},\bar{c}\text{ arbitrary},\ D\in C_{\mathrm{GL}_3(K)}(W)\right\},&x\in\reviewtext{pS_{m+1}}.
\end{cases}
\]
\begin{blueblock}
This follows directly from the full matching equation; when $x$ is a unit it forces $\bar{b}^{\,t}=\bar{c}=0$.

\end{blueblock}
\subsection{Counts for every \texorpdfstring{$m\geq1$}{m >= 1}}
\begin{blueblock}
For $|K|=q$, write $\newtext{\nu_2}(B)$ for the number of ordinary two-dimensional similarity classes of the indicated type and $h_2(B)$ for their common centralizer order. For scalar $B$ these are $q$ and $|\mathrm{GL}_2(K)|$; for companion $B$ they are $\nu_\theta,h_\theta$ from \eqref{eq:quadratic-data}. Similarly, the complete three-dimensional data are
{\color{black}
\begin{center}\small
\begin{tabular}{l c c}
\toprule
Type of $W$ & $\newtext{\nu_3}(W)$ & $h_3(W)$\\\midrule
Scalar & $q$ & $|\mathrm{GL}_3(K)|$\\
One eigenvalue, Jordan type $(2,1)$ & $q$ & $q^3(q-1)^2$\\
One eigenvalue, Jordan type $(3)$ & $q$ & $q^2(q-1)$\\
Two eigenvalues, semisimple & $q(q-1)$ & $q(q-1)^3(q+1)$\\
Two eigenvalues, nonsemisimple & $q(q-1)$ & $q(q-1)^2$\\
Three distinct eigenvalues & $q(q-1)(q-2)/6$ & $(q-1)^3$\\
Irreducible quadratic and a linear factor & $q^2(q-1)/2$ & $(q^2-1)(q-1)$\\
Irreducible cubic & $(q^3-q)/3$ & $q^3-1$\\\bottomrule
\end{tabular}
\end{center}
}
The parameter counts follow by choosing the eigenvalues and monic irreducible polynomials. The centralizer orders follow from the corresponding $K[T]$-endomorphism algebras, or from \eqref{eq:field-centralizer} below. In particular,
\[
 \sum\newtext{\nu_2}=q^2+q,\quad\sum\newtext{\nu_3}=q^3+q^2+q,\qquad
 \sum\frac{\newtext{\nu_2}}{h_2}=\frac{q^4}{|\mathrm{GL}_2(K)|},\quad
 \sum\frac{\newtext{\nu_3}}{h_3}=\frac{q^9}{|\mathrm{GL}_3(K)|}.
\]
\end{blueblock}
\edittext{Put $G_m=|G_\lambda|=q^{m+9}(q-1)^2(q^2-1)(q^3-1)$. We use the table conventions of Section~\ref{subsec:extension-action}. Each row involving $B$ or $W$ contributes once for each ordinary similarity type specified above. In rows marked ``each'', the entries apply separately to both listed types.}
\afterpage{\begin{landscapetablepages}
\begin{center}\normalfont\small
Orbit counts and the class equation for $\lambda=(m+1,1,1,1)$, $m\geq1$

\medskip
\renewcommand{\arraystretch}{1.6}
\begin{tabular}{c c c c c}
\toprule
Type & Number $N$ & Stabilizer order $H$ & $|\mathcal O|=G_m/H$ & $N|\mathcal O|$\\\midrule
$1,2\ (z\neq0)$ & $q^m(q-1)\newtext{\nu_2}$ & $\newtext{q^{m+1}}(q-1)h_2$ & $\dfrac{q^8(q-1)(q^2-1)(q^3-1)}{h_2}$ & $\dfrac{\newtext{q^{m+8}}(q-1)^2(q^2-1)(q^3-1)\newtext{\nu_2}}{h_2}$\\[8pt]
$3$ & $\newtext{q^{m+1}}$ & $\newtext{q^{m+3}}(q-1)$ & $q^6(q-1)(q^2-1)(q^3-1)$ & $\newtext{q^{m+7}}(q-1)(q^2-1)(q^3-1)$\\
$4$ & $q^m$ & $\newtext{q^{m+3}}(q-1)$ & $q^6(q-1)(q^2-1)(q^3-1)$ & $\newtext{q^{m+6}}(q-1)(q^2-1)(q^3-1)$\\
$5$ & $\newtext{q^{m+1}}$ & $\newtext{q^{m+4}}(q-1)^2$ & $q^5(q^2-1)(q^3-1)$ & $\newtext{q^{m+6}}(q^2-1)(q^3-1)$\\
$6$ & $\newtext{q^{m+1}}(q-1)$ & $\newtext{q^{m+2}}(q-1)^2$ & $q^7(q^2-1)(q^3-1)$ & $\newtext{q^{m+8}}(q-1)(q^2-1)(q^3-1)$\\
$1,2\ (z=0)$ & $q^m\newtext{\nu_2}$ & $\newtext{q^{m+5}}(q-1)h_2$ & $\dfrac{q^4(q-1)(q^2-1)(q^3-1)}{h_2}$ & $\dfrac{\newtext{q^{m+4}}(q-1)(q^2-1)(q^3-1)\newtext{\nu_2}}{h_2}$\\[8pt]
$9,10$ & $q^m\newtext{\nu_2}$ & $\newtext{q^{m+5}}(q-1)h_2$ & $\dfrac{q^4(q-1)(q^2-1)(q^3-1)}{h_2}$ & $\dfrac{\newtext{q^{m+4}}(q-1)(q^2-1)(q^3-1)\newtext{\nu_2}}{h_2}$\\[8pt]
$7,11$ (each) & $q^m$ & $\newtext{q^{m+3}}(q-1)$ & $q^6(q-1)(q^2-1)(q^3-1)$ & $\newtext{q^{m+6}}(q-1)(q^2-1)(q^3-1)$\\
$8,12$ (each) & $\newtext{q^{m+1}}$ & $\newtext{q^{m+4}}(q-1)^2$ & $q^5(q^2-1)(q^3-1)$ & $\newtext{q^{m+6}}(q^2-1)(q^3-1)$\\
$13,\ x\in \reviewtext{S_{m+1}^\times}$ & $q^m(q-1)\newtext{\nu_3}$ & $q^m(q-1)h_3$ & $\dfrac{q^9(q-1)(q^2-1)(q^3-1)}{h_3}$ & $\dfrac{\newtext{q^{m+9}}(q-1)^2(q^2-1)(q^3-1)\newtext{\nu_3}}{h_3}$\\[8pt]
$13,\ x\in\reviewtext{pS_{m+1}}$ & $q^m\newtext{\nu_3}$ & $\newtext{q^{m+6}}(q-1)h_3$ & $\dfrac{q^3(q-1)(q^2-1)(q^3-1)}{h_3}$ & $\dfrac{\newtext{q^{m+3}}(q-1)(q^2-1)(q^3-1)\newtext{\nu_3}}{h_3}$\\[8pt]\bottomrule
\end{tabular}
\end{center}
\end{landscapetablepages}}\par
\clearpage
\begin{blueblock}
\edittext{The stabilizer orders are obtained by counting the free parameters above. The $m$ higher coefficients of $a$ give $q^m$ choices; the remaining conditions are over $K$. Compared with $m=1$, each representative and each stabilizer has $m-1$ additional free coefficients, in $x$ and $a$, respectively. Thus both $N$ and $H$ are multiplied by $q^{m-1}$, while the orbit sizes stay the same.}

Summing the product column gives $\sum N G_m/H=\newtext{q^{m+16}}=|E_{(m+1,1,1,1)}|$, while summing the representative counts gives
\begin{equation}\label{eq:chook}
 c_{(m+1,1,1,1)}(q)=q^m(q^4+2q^3+4q^2+4q+3).
\end{equation}
The partition $(2,1,1,1)$ is obtained by taking $m=1$.

\end{blueblock}
\begin{addedblock}
\section{The classification in dimension five}\label{sec:completion}
\subsection{The known primary cases}\label{subsec:known-primary}
\redtext{The cyclic and elementary cases are given by \cite[Theorems 8.1 and 8.2]{PSS}. For $\lambda=(r)$, the $G_\lambda(p)$-action on $E_\lambda(p)=S_r$ is trivial. For $\lambda=(1^r)$, the action is the usual conjugation action of $\mathrm{GL}_r(K)$ on $M_r(K)$.}

In the elementary case, choose the usual rational canonical representatives over $K$. More explicitly, for a new polynomial variable $u$, write the corresponding $K[u]$-module as
\[
 \bigoplus_f\bigoplus_i K[u]/(f^{\mu_i(f)}),\qquad
 \sum_f\deg(f)|\mu(f)|=r.
\]
Its automorphism group is the stabilizer. \redtext{If $|K|=q$, its order is given by \cite[Theorem 9.5]{PSS}:}
\begin{equation}\label{eq:field-centralizer}
 \prod_f\left[ q^{\deg(f)\sum_j(\mu_j'(f))^2}
 \prod_{a\geq1}\prod_{j=1}^{m_a(\mu(f))}
       (1-q^{-j\deg(f)})\right],
\end{equation}
where $\mu'$ is the conjugate partition and $m_a(\mu)$ is the multiplicity of the part $a$.

For $|\lambda|=5$, the sources of the required primary data are as follows.
\begin{center}
\begin{tabular}{lll}
\toprule
Partition & Orbit classification & $\sum_{i,j}\min(\lambda_i,\lambda_j)$\\
\midrule
$(5)$ & \redtext{\cite[Theorem 8.1]{PSS}} & $5$\\
$(4,1)$ & \cite[Theorem 8.5]{PSS}, $m=3$ & $7$\\
$(3,2)$ & Section~\ref{sec:32} & $9$\\
$(3,1,1)$ & \cite[Theorem 8.6]{PSS}, $m=2$ & $11$\\
$(2,2,1)$ & Section~\ref{sec:221} & $13$\\
$(2,1,1,1)$ & Section~\ref{sec:hook}, $m=1$ & $17$\\
$(1^5)$ & \redtext{\cite[Theorem 8.2]{PSS}} & $25$\\
\bottomrule
\end{tabular}
\end{center}
The last column supplies the exponent of $Q$ in the centralizer kernel when $p$ is linear. In particular, the stabilizer orders in the tables for $(3,2)$, $(2,2,1)$ and $(2,1,1,1)$ must be multiplied by $Q^9$, $Q^{13}$ and $Q^{17}$, respectively, to obtain the orders of the corresponding centralizers in $\mathrm{GL}_5(R)$.

For completeness, the smaller nonempty partitions that may occur as primary components are
\[
\begin{array}{c|l}
1&(1)\\
2&(2),(1,1)\\
3&(3),(2,1),(1,1,1)\\
4&(4),(3,1),(2,2),(2,1,1),(1,1,1,1).
\end{array}
\]
Besides the cyclic and elementary cases, these are covered by \cite[Theorems 8.3, 8.5 and 8.6]{PSS}; the rectangular partition $(2,2)$ reduces to similarity of two by two matrices over $S_2$, whose normal forms appear in Theorem~\ref{nf} in Section~\ref{sec:221}. For $\deg p\geq2$ and $\deg(p)|\lambda|\leq5$, necessarily $|\lambda|\leq2$. These higher-degree primary cases are therefore cyclic or elementary.

\subsection{\reviewtext{Completion of the classification}}\label{subsec:completion}
We prove the first part of Theorem~\ref{thm:main} and describe the corresponding centralizers. Choose the rational canonical class of the reduction $A\in M_5(k)$. It is specified by finitely many partitions $\lambda(p)$ satisfying $\sum_p\deg(p)|\lambda(p)|=5$. The lifting correspondence and primary decomposition reduce the choice of a similarity class above $A$ to an independent choice of one orbit in each $G_{\lambda(p)}(p)\backslash E_{\lambda(p)}(p)$.

If $|\lambda(p)|=5$, then $\deg p=1$ and the preceding table lists all seven possibilities. If $|\lambda(p)|<5$, the smaller-partition list and the cited results apply. Thus the three orbit classifications of this paper, together with those previous results, cover every primary component.

We use the construction of \cite[Section 6.4 and Theorem 6.5]{PSS} to turn each primary extension representative $E$ into a matrix over $R$. Let $d=\deg p$ and use the chosen monic lift $\widetilde p$ of $p$. For each entry $E_{ij}$, lift the coefficients of its polynomial representative of degree less than $d\min(\lambda_i,\lambda_j)$ to obtain $\widetilde E_{ij}$. Form the $(i,j)$ block of $X_E$, of size $d\lambda_i\times d\lambda_j$, by placing the coefficients of $\widetilde E_{ij}$ in increasing order of degree in the last column, padded with zeros. The other columns are zero. The matrix corresponding to $E$ is
\begin{equation}\label{eq:matrix-representative}
 \mathcal A_E=A_0+\pi X_E,\qquad
 A_0=\bigoplus_i\mathsf C(\widetilde p^{\lambda_i}),
\end{equation}
where $\mathsf C(f)$ \proofreadchange{for monic $f$} is the matrix of multiplication by $t$ on $R[t]/(f)$ in the basis $1,t,\ldots,t^{\deg f-1}$. The matrix $\pi X_E$ is independent of the coefficient lifts. For instance, when $p=t-a$ and $\lambda=(3,2)$, only columns $3$ and $5$ of $\pi X_E$ can be nonzero.

Apply this construction to each primary component and take their direct sum, with the primary blocks in a fixed order. This produces a five by five matrix over $R$ in every similarity class. Two such matrices can be similar only if their reductions have the same primary data and their extension matrices belong to the same primary orbits. The comparisons of representatives in Sections~\ref{sec:32}--\ref{sec:hook}, together with the known primary classifications, then show that the parameter choices agree.

The centralizers follow from \eqref{eq:centralizer-primary-product} and \eqref{eq:centralizer-sequence}. For a primary block $\mathcal A_E\in M_n(R)$, put $A=\rho(\mathcal A_E)$. For $g\in\Stab_{G_A}(E)$, viewed as an $n\times n$ matrix over $k$, choose a lift $\widetilde g\in M_n(R)$ and define $D_g\in M_n(k)$ by $\mathcal A_E\widetilde g-\widetilde g\mathcal A_E=\pi\widetilde D_g$. The centralizer elements reducing to $g$ are
\begin{equation}\label{eq:lift-centralizer}
 \widetilde g+\pi\widetilde Y,\qquad AY-YA=-D_g.
\end{equation}
The exact sequence \eqref{eq:centralizer-sequence} ensures that this equation has a solution for every $g$ in the stabilizer. The finite-field count is proved in the next section.

\section{Class sizes and the number of similarity classes}\label{sec:enumeration}
Assume that $|k|=Q$. We first collect the primary orbit counts and then combine them using primary decomposition.

\subsection{Primary orbit counts}\label{subsec:primary-counts}
\edittext{The primary orbit counts are listed below. The three new entries are \eqref{eq:c32}, \eqref{eq:c221} and \eqref{eq:chook}, with $m=1$ in the last formula. The other entries come from the primary classifications recalled in Section~\ref{sec:completion}.}

\begin{center}\small
\begin{tabular}{c l}
\toprule
$\lambda$ & $c_\lambda(q)$\\\midrule
$(1)$ & $q$\\
$(2)$ & $q^2$\\
$(1,1)$ & $q^2+q$\\
$(3)$ & $q^3$\\
$(2,1)$ & $q^3+q^2+q$\\
$(1,1,1)$ & $q^3+q^2+q$\\
$(4)$ & $q^4$\\
$(3,1)$ & $q^4+q^3+q^2$\\
$(2,2)$ & $q^4+q^3+q^2$\\
$(2,1,1)$ & $q^4+2q^3+2q^2+2q$\\
$(1^4)$ & $q^4+q^3+2q^2+q$\\
$(5)$ & $q^5$\\
$(4,1)$ & $q^5+q^4+q^3$\\
$(3,2)$ & $q^5+q^4+2q^3+q^2$\\
$(3,1,1)$ & $q^5+2q^4+2q^3+2q^2$\\
$(2,2,1)$ & $q^5+2q^4+3q^3+2q^2+2q$\\
$(2,1,1,1)$ & $q^5+2q^4+4q^3+4q^2+3q$\\
$(1^5)$ & $q^5+q^4+2q^3+2q^2+q$\\\bottomrule
\end{tabular}
\end{center}
\edittext{The cyclic, elementary and rectangular entries follow from the results recalled in Section~\ref{sec:completion}. The entries for $(2,1)$, $(3,1)$, $(4,1)$, $(2,1,1)$ and $(3,1,1)$ come from \cite[Tables 8.1 and 8.2]{PSS}.}

\subsection{Combining the primary components}
The number of monic irreducible polynomials of degree $d$ over $\mathbb F_Q$ is
\[
 I_d(Q)=\frac1d\sum_{e\mid d}\mu(e)Q^{d/e},
\]
where $\mu$ is the M\"obius function. Write $[z^j]F(z)$ for the coefficient of $z^j$ in a formal power series $F(z)$.

\redtext{Specializing the generating-function formula in \cite[Section 9, p.~510]{PSS} to the coefficient of degree five gives}
\begin{equation}\label{eq:count-product}
 s_5(Q)=[z^5]\prod_{d=1}^{5}
 \left(\sum_{|\lambda|\leq\lfloor5/d\rfloor}
 c_\lambda(Q^d)z^{d|\lambda|}\right)^{I_d(Q)}.
\end{equation}
\redtext{Here $s_5(Q)$ is the number of similarity classes in $M_5(R)$. We evaluate this formula below using the primary counts in Section~\ref{subsec:primary-counts}.}

\edittext{For a representative $\mathcal A$ with primary data $(\lambda(p),E_p)$, \proofreadchange{put $A=\rho(\mathcal A)$.} \proofreadchange{Equation}~\eqref{eq:global-centralizer} and the orbit--stabilizer formula give}
\begin{equation}\label{eq:class-size}
 |\mathrm{GL}_5(R)\cdot\mathcal A|
 =\frac{Q^{25}|\mathrm{GL}_5(\mathbb F_Q)|}
 {Q^{\delta(A)}\prod_p|\Stab_{G_{\lambda(p)}(p)}(E_p)|},
 \qquad |\mathrm{GL}_5(\mathbb F_Q)|=\prod_{j=0}^{4}(Q^5-Q^j).
\end{equation}

\begin{cor}\label{cor:total-five}
For every finite commutative local principal ideal ring $R$ of length two with residue field of order $Q$, the number of similarity classes in $M_5(R)$ is
\begin{equation}\label{eq:total-five}
 s_5(Q)=Q^{10}+Q^9+3Q^8+4Q^7+8Q^6+8Q^5+10Q^4+6Q^3+6Q^2.
\end{equation}
In particular, this number depends only on the residue-field order.
\end{cor}
\begin{proof}
Put $A_j(q)=\sum_{|\lambda|=j}c_\lambda(q)$, with $A_0(q)=1$. Summing the rows of the table gives
\begin{align*}
 A_1(q)&=q,\\
 A_2(q)&=2q^2+q,\\
 A_3(q)&=3q^3+2q^2+2q,\\
 A_4(q)&=5q^4+5q^3+6q^2+3q,\\
 A_5(q)&=7q^5+9q^4+14q^3+11q^2+6q.
\end{align*}
Let $F(z)=(1+\sum_{j=1}^5 A_j(Q)z^j)^Q$ be the contribution from the $Q$ monic linear polynomials, and put $B_d=I_d(Q)Q^d$ for $2\le d\le5$. For $d=3,4,5$, only partitions of size one give nonconstant terms of degree at most five. Since $A_1(Q^d)=Q^d$, the factor for degree $d$ in \eqref{eq:count-product} is $1+B_dz^d\proofreadchange{\pmod{z^6}}$.

For the quadratic polynomials, the term in $z^4$ comes from either one irreducible quadratic with a partition of size two or two distinct irreducible quadratics with a partition of size one for each. In the second case, each quadratic contributes $A_1(Q^2)=Q^2$ choices. Thus the factor for degree two is $1+B_2z^2+T_4z^4\proofreadchange{\pmod{z^6}}$, where
\[
 T_4=I_2(Q)A_2(Q^2)+\binom{I_2(Q)}{2}Q^4.
\]
Multiplying the factors for degrees two through five gives
\[
 1+B_2z^2+B_3z^3+(T_4+B_4)z^4+(B_5+B_2B_3)z^5
 \pmod{z^6}.
\]
Taking the coefficient of $z^5$ after multiplying by $F(z)$ gives
\[
 s_5(Q)=[z^5]F+B_2[z^3]F+B_3[z^2]F
       +(T_4+B_4)[z]F+B_5+B_2B_3.
\]
\Needspace{6\baselineskip}
Write $a_j=A_j(Q)$ for $1\le j\le5$. To compute $[z^j]F$, we choose one term from each of the $Q$ factors of $F(z)$ so that their degrees sum to $j$. For example, in degree five the choice of degrees $3,1,1$ contributes $3\binom{Q}{3}a_3a_1^2$, since any of the three chosen factors can supply the term of degree three. The required coefficients are
\begin{align*}
 [z]F &= Qa_1,\qquad [z^2]F=Qa_2+\binom{Q}{2}a_1^2,\\
 [z^3]F &= Qa_3+2\binom{Q}{2}a_2a_1+\binom{Q}{3}a_1^3,\\
 [z^5]F &= Qa_5+2\binom{Q}{2}(a_4a_1+a_3a_2)
           +3\binom{Q}{3}(a_3a_1^2+a_2^2a_1)\\*
        &\qquad {}+4\binom{Q}{4}a_2a_1^3+\binom{Q}{5}a_1^5.
\end{align*}
Substituting the values of $A_j$ and $I_d(Q)$ into the coefficient formula gives \eqref{eq:total-five}.
\end{proof}

\section*{Acknowledgements}
\textbf{The research of Gunja Sachdeva and the Junior Research Fellowship of Adeep Varghese are supported by the Anusandhan National Research Foundation (ANRF), a statutory body of the Department of Science and Technology, Government of India, under SERB POWER Grant No.~SPG/2022/001738.}

\smallskip
\noindent\begin{minipage}[t]{\linewidth}
\footnotesize\color{black}
\textsc{Department of Mathematics, BITS Pilani, K K Birla Goa Campus}\\
\textsc{Zuarinagar, Goa 403726, India}\par
\smallskip
\begin{tabular}{@{}l@{\quad}l@{}}
Himadri Mukherjee &
  \href{mailto:himadrim@goa.bits-pilani.ac.in}{\texttt{himadrim@goa.bits-pilani.ac.in}}\\
Gunja Sachdeva &
  \href{mailto:gunjas@goa.bits-pilani.ac.in}{\texttt{gunjas@goa.bits-pilani.ac.in}}\\
Adeep Varghese &
  \href{mailto:p20230031@goa.bits-pilani.ac.in}{\texttt{p20230031@goa.bits-pilani.ac.in}}
\end{tabular}
\end{minipage}\par

\end{addedblock}

\end{document}